
\def\LaTeX{%
  \let\Begin\begin
  \let\End\end
  \let\salta\relax
  \let\finqui\relax
  \let\futuro\relax}

\def\UK{\def\our{our}\let\sz s}
\def\USA{\def\our{or}\let\sz z}

\UK



\LaTeX

\USA


\salta

\documentclass[twoside,12pt]{article}
\setlength{\textheight}{24cm}
\setlength{\textwidth}{16cm}
\setlength{\oddsidemargin}{2mm}
\setlength{\evensidemargin}{2mm}
\setlength{\topmargin}{-15mm}
\parskip2mm


\usepackage[usenames,dvipsnames]{color}
\usepackage{amsmath}
\usepackage{amsthm}
\usepackage{amssymb}
\usepackage[mathcal]{euscript}


\usepackage[T1]{fontenc}
\usepackage[latin1]{inputenc}
\usepackage[english]{babel}
\usepackage[babel]{csquotes}

\usepackage{cite}

\usepackage{latexsym}
\usepackage{graphicx}
\usepackage{mathrsfs}
\usepackage{mathrsfs}
\usepackage{hyperref}
\usepackage{pgfplots}

%
%


\definecolor{viola}{rgb}{0.3,0,0.7}
\definecolor{ciclamino}{rgb}{0.5,0,0.5}

\def\gianni #1{{\color{blue}#1}}
\def\pier #1{{\color{red}#1}}
\def\juerg #1{{\color{green}#1}}
\def\juergen #1{{\color{red}#1}}

\def\gianni #1{#1}
\def\pier #1{#1}
\def\juerg #1{#1}
\def\juergen #1{#1}




\bibliographystyle{plain}


%

\finqui

\def\Beq{\Begin{equation}}
\def\Eeq{\End{equation}}
\def\Bsist{\Begin{eqnarray}}
\def\Esist{\End{eqnarray}}

\def\Bthm{\Begin{theorem}}
\def\Ethm{\End{theorem}}
\def\Blem{\Begin{lemma}}
\def\Elem{\End{lemma}}

\def\Bcor{\Begin{corollary}}
\def\Ecor{\End{corollary}}
\def\Brem{\Begin{remark}\rm}
\def\Erem{\End{remark}}

\def\Bdim{\Begin{proof}}
\def\Edim{\End{proof}}
\def\Bcenter{\Begin{center}}
\def\Ecenter{\End{center}}
\let\non\nonumber




\def\step #1 \par{\medskip\noindent{\bf #1.}\quad}


\def\aand{\quad\hbox{and}\quad}

\def\lhs{left-hand side}
\def\rhs{right-hand side}



\def\multibold #1{\def\arg{#1}%
  \ifx\arg\pto \let\next\relax
  \else
  \def\next{\expandafter
    \def\csname #1#1#1\endcsname{{\bf #1}}%
    \multibold}%
  \fi \next}

\def\pto{.}

\def\multical #1{\def\arg{#1}%
  \ifx\arg\pto \let\next\relax
  \else
  \def\next{\expandafter
    \def\csname cal#1\endcsname{{\cal #1}}%
    \multical}%
  \fi \next}


\def\multimathop #1 {\def\arg{#1}%
  \ifx\arg\pto \let\next\relax
  \else
  \def\next{\expandafter
    \def\csname #1\endcsname{\mathop{\rm #1}\nolimits}%
    \multimathop}%
  \fi \next}

\multibold
qwertyuiopasdfghjklzxcvbnmQWERTYUIOPASDFGHJKLZXCVBNM.

\multical
QWERTYUIOPASDFGHJKLZXCVBNM.

\multimathop
diag dist div dom mean meas sign supp .


\def\accorpa #1#2{\eqref{#1}--\eqref{#2}}
\def\Accorpa #1#2 #3 {\gdef #1{\eqref{#2}--\eqref{#3}}%
  \wlog{}\wlog{\string #1 -> #2 - #3}\wlog{}}


\def\separa{\noalign{\allowbreak}}

\def\somma #1#2#3{\sum_{#1=#2}^{#3}}

\def\<#1>{\mathopen\langle #1\mathclose\rangle}
\def\norma #1{\mathopen \| #1\mathclose \|}

\def\[#1]{\mathopen\langle\!\langle #1\mathclose\rangle\!\rangle}

\def\iot {{\int_0^t}}
\def\ioT {{\int_0^T}}
\def\intQt{{\int_0^t\!\!\int_\Omega}}

\def\iO{\int_\Omega}

\def\dt{\partial_t}
\def\dn{\partial_{\bf n}}

\def\cpto{\,\cdot\,}

\def\checkmmode #1{\relax\ifmmode\hbox{#1}\else{#1}\fi}

\def\aet{\checkmmode{a.e.\ in~$(0,T)$}}

\def\aat{\checkmmode{for a.a.~$t\in(0,T)$}}


\def\erre{{\mathbb{R}}}

\def\enne{{\mathbb{N}}}




\def\genspazio #1#2#3#4#5{#1^{#2}(#5,#4;#3)}
\def\spazio #1#2#3{\genspazio {#1}{#2}{#3}T0}
\def\spaziot #1#2#3{\genspazio {#1}{#2}{#3}t0}

\def\L {\spazio L}
\def\H {\spazio H}
\def\W {\spazio W}
\def\Lt {\spaziot L}

\def\C #1#2{C^{#1}([0,T];#2)}


\def\Lx #1{L^{#1}(\Omega)}
\def\Hx #1{H^{#1}(\Omega)}

\def\LQ #1{L^{#1}(Q)}

\def\Luno{\Lx 1}
\def\Ldue{\Lx 2}
\def\Linfty{\Lx\infty}

\def\Huno{\Hx 1}
\def\Hdue{\Hx 2}
\def\Hunoz{{H^1_0(\Omega)}}


\def\LQ #1{L^{#1}(Q)}

\def\Liq{L^\infty(Q)}

\let\theta\vartheta
\let\eps\varepsilon
\let\phi\varphi

\let\hat\widehat

\let\TeXchi\chi                         
\newbox\chibox
\setbox0 \hbox{\mathsurround0pt $\TeXchi$}
\setbox\chibox \hbox{\raise\dp0 \box 0 }
\def\chi{\copy\chibox}


\def\pn{\mathbf{n}}

\def\VA #1{V_A^{#1}}
\def\VB #1{V_B^{#1}}

\def\F1l{F_{1,\lambda}}

\def\phil{\phi_\lambda}
\def\thetal{\theta_\lambda}

\def\qn{q_n}
\def\pen{p_n}

\def\ualn{u_{\alpha_n}}
\def\tal{\theta_\alpha}
\def\pal{\phi_\alpha}
\def\hal{h_\alpha}
\def\dhal{h'_\alpha}

\def\taln{\theta_{\alpha_n}}
\def\paln{\phi_{\alpha_n}}

\def\haln{h_{\alpha_n}}
\def\phih{\phi^h}
\def\thetah{\theta^h}
\def\etah{\eta^h}
\def\xih{\xi^h}
\def\bphi{\overline\phi}
\def\bu{\overline u}
\def\btheta{\overline \theta}
\def\etan{\eta_n}
\def\xin{\xi_n}
\def\un{u_n}
\def\thetan{\theta_n}
\def\phin{\phi_n}
\def\uad{{\cal U}_{\rm ad}}
\def\J{{\cal J}}
\def\indi{I_{[-1,1]}}
\Begin{document}


%
\title{Optimal control of a phase field system of\\ Caginalp type with fractional operators}
\author{}
\date{}

\maketitle
\Bcenter
\vskip-1cm
{\large\sc Pierluigi Colli$^{(1)}$}\\
{\normalsize e-mail: {\tt pierluigi.colli@unipv.it}}\\[.25cm]
{\large\sc Gianni Gilardi$^{(1)}$}\\
{\normalsize e-mail: {\tt gianni.gilardi@unipv.it}}\\[.25cm]
{\large\sc J\"urgen Sprekels$^{(2)}$}\\
{\normalsize e-mail: {\tt sprekels@wias-berlin.de}}\\[.45cm]
$^{(1)}$
\pier{\small Dipartimento di Matematica ``F. Casorati''}\\
\pier{\small Universit\`a di Pavia}\\
\pier{\small via Ferrata 5, 27100 Pavia, Italy}\\[.45cm]
$^{(2)}$
{\small Department of Mathematics}\\
{\small Humboldt-Universit\"at zu Berlin}\\
{\small Unter den Linden 6, 10099 Berlin, Germany}\\[2mm]
{\small and}\\[2mm]
{\small Weierstrass Institute for Applied Analysis and Stochastics}\\
{\small Mohrenstrasse 39, 10117 Berlin, Germany}
\Ecenter
%
%
%
\Begin{abstract}\noindent
In their recent work ``Well-posedness, regularity and asymptotic analyses for a fractional phase field system''  
(\emph{Asymptot.\ Anal.} {\bf 114} (2019), 93--128), two of the present authors have studied phase field systems of Caginalp type, which model nonconserved, nonisothermal phase transitions and  in which the occurring diffusional operators are given by fractional versions in the spectral sense of unbounded, monotone, selfadjoint, linear operators having compact resolvents.  
In this paper, we complement this analysis by investigating distributed optimal control problems for such systems. 
It is shown that the associated control-to-state operator is Fr\'echet differentiable between suitable Banach spaces, 
and meaningful first-order necessary optimality conditions are derived in terms of a variational inequality and the associated adjoint state variables. 

 
\End{abstract}
\salta
\pagestyle{myheadings}
\newcommand\testopari{\sc Colli \ --- \ Gilardi \ --- \ Sprekels}
\newcommand\testodispari{\sc Optimal control of fractional Caginalp systems}
\markboth{\testopari}{\testodispari}
\finqui
%

\section{Introduction}
\label{Intro}
\setcounter{equation}{0}
The \emph{Caginalp phase field model} is a well-known system of partial differential equations modeling
the evolution of a temperature-dependent phase transition with nonconserved order parameter
$\,\phi\,$ that takes place in a container 
$\Omega\subset \erre^3$. A classical form that was introduced and analyzed in the seminal paper \cite{Cp1} 
(see also, e.g., \cite{Cp2,Cp3,Cp4}) is given by the evolutionary system 
\begin{align*}
\rho\,C_V\,\dt\theta\,+\,\ell \,\dt\varphi\,-\,\kappa\,\Delta\theta\,&=\,u,\\[1mm]
\alpha\,\xi^2\,\dt\varphi\,-\,\xi^2\,\Delta\varphi\,+\,F'(\phi)\,&=\,2\,\theta,
\end{align*} 
which is to be satisfied in the set $Q:=\Omega\times (0,T)$, where $T>0$ is a given final time.  
The first equation in the above system is an approximation to the universal 
balance law of internal energy, while the second one 
governs the evolution of the order parameter. 
The quantities $\rho,C_V,\ell,\kappa,\alpha,\xi$ are positive physical
constants; in particular, $\ell$ is closely allied to the latent heat released or absorbed during the phase transition 
process, and $\,\xi\,$ is a measure for the thickness of the transition zone between the different phases. 
The unknowns $\theta$ and $\phi$ stand for a temperature difference and the order
parameter (usually a normalized fraction of one of the phases involved in the phase transition), while
$u$ represents a control (a heat source or sink) and $\gianni F$ is a double-well potential 
whose derivative $\gianni{F'}$ is the thermodynamic force driving the phase transition. 
For a derivation of the model equations from general thermodynamic principles, we refer the reader to 
\cite[Chapter 4]{BS}.    

In their recent paper \cite{CG}, two of the present authors have studied a variation of the Caginalp 
model, namely the system 
\begin{align}
\label{sstheta}
&\partial_t\theta\,+\,\ell(\varphi)\partial_t\varphi\,+\,A^{2\rho}\theta\,=\,u\quad\mbox{in \,$Q$,}\\[0.5mm]
\label{ssphi}
&\partial_t\varphi \,+\,B^{2\sigma}\varphi\,+\,F'(\varphi)\,=\,\theta\,\ell(\varphi) \quad\mbox{in \,$Q$,}\\[0.5mm]
\label{ssini}
&\theta(0)\,=\,\theta_0,\quad \varphi(0)\,=\,\varphi_0,\quad\mbox{in \,$\Omega$}.
\end{align}
The main difference to the Caginalp system (besides the fact that many physical constants are 
normalized to unity and that the quantity $\ell$ representing the latent heat is allowed to depend
on the order parameter~$\phi$) is given by the fact that in \eqref{sstheta}--\eqref{ssphi}
the expressions  $A^{2\rho}$ and~$B^{2\sigma}$, with $\rho>0$ and $\sigma>0$,
denote fractional powers in the spectral sense of self-adjoint,
monotone, and unbounded linear operators $A$ and~$B$, respectively, 
which are supposed to be densely defined in $H:=\Ldue$ and to have compact resolvents. 
The standard example occurs when $A^{2\rho}=B^{2\sigma}=-\Delta$, with zero Dirichlet or Neumann boundary conditions.

The nonlinearity $\ell\,$ is assumed to be a smooth function, while\, $F$~denotes a double-well potential.
Typical and physically significant examples are the so-called {\em classical regular potential}, 
the {\em logarithmic potential\/},
and the {\em double obstacle potential\/}, which are given, in this order,~by
\begin{align}
\label{regpot}
  & F_{\rm reg}(r) := \frac 14 \, (r^2-1)^2 \,,
  \quad r \in \erre,\\[1mm] 
  \separa
  \label{logpot}
  & F_{\rm log}(r) := \left\{\begin{array}{ll}
   (1+r)\ln (1+r)+(1-r)\ln (1-r) - c_1 r^2\,,&
  \quad r \in (-1,1)\\
  2\log(2)-c_1\,,&\quad r\in\{-1,1\}\\
  +\infty\,,&\quad r\not\in [-1,1]
\end{array}\right. \,,
  \\[1mm]
  \separa
\label{obspot}
  & F_{\rm 2obs}(r) :=  c_2(1- r^2) 
  \quad \hbox{if $|r|\leq1$},
  \aand
  F_{\rm 2obs}(r) := +\infty
  \quad \hbox{if $|r|>1$}.
\end{align}
Here, the constants  in \eqref{logpot} and \eqref{obspot} satisfy
$c_1>1$ and $c_2>0$, so that the corresponding functions are nonconvex.
In cases like \eqref{obspot}, one has to split $F$ into a nondifferentiable convex part~$F_1$ 
(the~indicator function of $[-1,1]$, in the present example) and a smooth perturbation~$F_2$.
Accordingly, in the term $F'(\phi)$ appearing in~\eqref{ssphi}, 
one has to replace the derivative $F_1'$ of the convex part $F_1$
by the subdifferential $\partial F_1$ and interpret \eqref{ssphi} as a differential inclusion
or as a variational inequality involving $F_1$ rather than~$\partial F_1$.

In \cite{CG}, general results on well-posedness, regularity and asymptotic behavior have been
proved for the state system \eqref{sstheta}--\eqref{ssini}. 
In this paper, we complement the analysis in \cite{CG}
by studying the optimal control of this system. 
More precisely, 
given nonnegative constants $\beta_i$, $1\le i\le 5$, target functions
$\varphi_\Omega,\theta_\Omega\in\Ldue$ and $\varphi_Q,\theta_Q\in L^2(Q)$, as well as
threshold functions $u_{min},u_{max}\in L^\infty(Q)$ with $\,u_{min}\le u_{max}\,$ in~$Q$, 
we consider the following optimal control problem:

\vspace{2mm}\noindent
{\bf (CP)}  \,\,\,Minimize the tracking-type cost functional
\begin{align}
\label{cost}
\J((\varphi,\theta),u)\,:=&\,\frac {\beta_1}2\iO|\varphi(T)-\varphi_\Omega|^2\,+\,\frac{\beta_2}2
\ioT\!\!\iO|\varphi-\varphi_Q|^2\,+\,\frac{\beta_3}2\iO|\theta(T)-\theta_\Omega|^2\nonumber\\
&\,+\,\frac{\beta_4}2 \ioT\!\!\iO|\theta-\theta_Q|^2\,+\,\frac{\beta_5}2\int_Q|u|^2
\end{align}
over the set of admissible controls
\begin{equation}
\label{uad}
\uad\,:=\,\left\{u\in L^ \infty(Q):\,u_{min}\,\le u\,\le\,u_{max} \,\mbox{ a.e.\ in }\,Q\right\},
\end{equation}
subject to the state system \eqref{sstheta}--\eqref{ssini}.

The optimal control problem {\bf (CP)} constitutes a generalization of investigations for the
original Caginalp phase field system with regular potential $F_{\rm reg}$ that were begun  
in the early nineties of the past century; in this
connection, we refer the reader to the pioneering works \cite{ChenHoff,Hein, HeinSachs,HeinTr,HoffJiang}
(see also the related sections in the monograph~\cite{Tr}). 
For more recent contributions, we mention the papers\pier{\cite{BCGMR,CGMR1,CGMR2,CM,LS}}, 
where in \cite{LS} a thermodynamically consistent version of the phase field system
was considered. 
We also mention the papers \cite{CMR,HSS,SZheng} \pier{that} were devoted to optimal control problems
for the Penrose--Fife phase field model of phase transitions with nonconserving kinetics.     

The problem {\bf (CP)} can also be seen in comparison with a class of optimal control
problems for Cahn--Hilliard type systems of the form 
\begin{align}\label{CH1}
&\alpha\dt\mu +\dt\phi+A^{2\rho}\mu =0,\\
\label{CH2}
&\beta\dt\phi+B^{2\sigma}\phi+F'(\phi)=\mu+u,\\
\label{CH3}
&\mu(0)=\mu_0,\quad\phi(0)=\phi_0,
\end{align} 
where $\mu$ represents the chemical potential and $\alpha\ge0$ and $\beta\ge0$. 
Obviously, 
\eqref{sstheta}--\eqref{ssini} is in the special case $\ell(\phi)\equiv \ell>0$ of the above type 
(put $\mu=\theta$, $\alpha=1/\ell$ and $\beta=1$), 
where, however, the control \,$u\,$ appears in the phase equation. 
For the case when $\alpha=0$ and $\beta>0$, optimal control problems 
for \eqref{CH1}--\eqref{CH3} have been treated in the recent papers \cite{CGS19, CGS21} \pier{and reviewed in \cite{CGS21bis}}.
Moreover, for the classical case when $A=B=-\Delta$, $\rho=\sigma=1/2$, with various boundary conditions 
(i.e., Dirichlet, Neumann, and dynamic conditions), 
there exist many contributions \gianni{in which} optimal control problems have been studied; 
for a number of recent references in this direction, we refer the reader to~\cite{CGS21}. 

Another closely related phase field system is given by a model for tumor growth for which control
problems have recently been studied. The fractional version of this model reads as follows
(cf.\ \cite{CGS25,CGSASY}):
\begin{align}\label{TG1}
&\alpha\dt\mu+\dt\phi+A^{2\rho}\mu\,=\,P(\phi)(\juergen{S}-\mu),\\
\label{TG2}
&\beta\dt\phi+B^{2\sigma}\phi+F'(\phi)=\mu ,\\
\label{TG3}
&\dt S+C^{2\tau}S=-P(\phi)(S-\mu),\\
\label{TG4}
&\mu(0)=\mu_0,\quad \phi(0)=\phi_0,\quad S(0)=S_0.
\end{align}
Indeed, if $P(\phi)\equiv 0$, then \eqref{TG1}, \eqref{TG2} decouple from \eqref{TG3} and 
attain the form \eqref{CH1}, \eqref{CH2} for $u=0$. 
Also for the system \eqref{TG1}--\eqref{TG4}
optimal control problems have been studied for the classical case $A=B=C=-\Delta$,
$\rho=\sigma=\tau=1/2$, with zero Neumann boundary conditions. 
In this connection, we refer to the works \cite{CGRS, S_a,S_b,S_DQ, S}. 
In~\cite{CSS1}, also terms \juerg{modeling} chemotaxis were incorporated in the model. 
Even more involved models have been studied in \cite{EK,EK_ADV,GLR}.  
 
In this paper, in which we study the state system \eqref{sstheta}--\eqref{ssini}, we have to 
focus on the interplay between the nonlinearity $F$ and embedding properties of the 
domains of the involved operators. 
Quite surprisingly,
it turns out that \gianni{in our} case, where \gianni{$\alpha=1/\ell>0$ if we consider~\eqref{CH1}},
the situation is more delicate than in the 
abovementioned works where $\alpha=0$. 
The reason for this is that a proper treatment of the
nonlinear term $F'(\phi)$ in the optimal control problem 
(in particular, the derivation of results concerning Fr\'echet differentiability) 
makes it necessary that $F'(\phi)\in L^\infty(Q)$. 
This means, at least in the case of the irregular potentials $F_{\rm log}$ and $F_{\rm 2obs}$,
that the values attained by the solution component $\phi$ must be uniformly separated from
the critical values (in this case~$\pm 1$). 
To show such a separation condition, however, it is
somehow needed that the right-hand side $\ell(\phi)\theta$ of \eqref{ssphi}
($\mu$, in the case of \eqref{CH2} or \eqref{TG2}) is bounded. 
In terms of the expected 
regularities, this condition is more restrictive in our situation. 
Indeed, if $V_A^\rho=D(A^\rho)$
denotes the domain of the fractional operator $A^\rho$, then it turns out that 
\gianni{it maximally holds $\theta\in L^\infty(0,T;\VA\rho)$ in the case of the system \eqref{sstheta}--\eqref{ssini}, 
which corresponds to $\mu\in\L\infty{\VA\rho}$ for the system \eqref{CH1}--\eqref{CH3} with $\alpha>0$,
while one can recover the \juerg{better} regularity $\mu\in L^\infty(0,T;V_A^{2\rho})$ if $\alpha=0$}. 

It turns out that an appropriate   separation property 
(the condition \juerg{{\bf (GB)}} below) holds true in certain cases for regular potentials and 
singular potentials of the logarithmic type. 
For such cases, the Fr\'echet differentiability
can be shown (see Section~\ref{FRECHET}), and first-order necessary optimality conditions can be
derived (see Section~\ref{NECOPT}). 
The last Section~\ref{DQ} brings the derivation of first-order necessary
conditions also for the case of the double obstacle potential~$F_{\rm 2obs}$. 
In this case, 
where a separation condition like {\bf (GB)} cannot be expected to hold and where we
do not have Fr\'echet differentiability, we apply the so-called {\em deep quench approximation}, 
taking advantage of the results established in Section~\ref{NECOPT} for logarithmic potentials. 

Throughout this paper, we denote for a given Banach space $(X, \|\,\cdot\,\|_X)$ by $\,X^*\,$ the
dual space of $X$ and by $\,\langle\,\cdot\,,\,\cdot\,\rangle_X\,$ the duality product between
$\,X^*\,$ and~$\,X$. 
We will also make frequent use of the elementary Young inequality
\gianni{%
\begin{equation}
\label{Young}
a\,b\,\le\,\,\delta a^2\,+\,\frac 1{4\delta} \, b^2
\quad\mbox{for all \ $a,b\in\erre$ \ and \ $\delta>0$}.
\end{equation}%
}%

Finally, we denote by $W^{s,p}(\Omega)$ for $s\ge 0$ and $p\in [1,+\infty]$ the fractional Sobolev--Slobodeckij spaces
defined in, e.g.,~\cite{NPV}. 
We put $H^s(\Omega):=W^{s,2}(\Omega)$ for $s\ge 0$ and notice that for $s\ge0$ and
in three dimensions of space  we have the continuous embeddings (cf., e.g., \cite[Thms.~6.7,~8.2]{NPV})
\begin{align}
\label{sobslo1}
& H^s(\Omega)\subset L^q(\Omega)\quad\mbox{for }\,2s<3 \,\mbox{ and }\,1\le q\le 6/(3-2s),\\[1mm]
\label{sobslo2}
& H^s(\Omega)\subset C^0(\overline\Omega)\quad\mbox{for }\,2s>3.
\end{align}  
Observe that the latter embedding is compact, while the former is compact only for $1\le q< 6/(3-2s)$. In particular,
we have 
\begin{equation}
\label{sobslo3}
H^{2s}(\Omega)\subset \Lx 4\,\,\mbox{ if }\,\,s\ge 3/8
\gianni{\aand}
H^{4s}(\Omega)\subset \Lx 6\,\,\mbox{ if }\,\, \gianni{s\ge 1/4} \,.
\end{equation}  

 
\section{Statement of the problem and and the state system}
\label{STATEMENT}
\setcounter{equation}{0}
 
In this section, we state precise assumptions and notations and present some results for
the state system \eqref{sstheta}--\eqref{ssini}. 
Throughout this paper,
$\Omega\subset\erre^3$ is a bounded and connected \gianni{open} set with smooth boundary $\Gamma:=\partial\Omega$ and volume~$|\Omega|$. 
We denote by $\pn$ the outward unit normal vector field and by $\partial_\pn$ the outward normal derivative. 
We~set
\Beq
  H := \Ldue
  \Eeq
and denote by $\norma\cpto$ and $(\cpto,\cpto)$ the standard norm and inner product of~$H$.
We generally assume:

\vspace{2mm}\noindent
{\bf (A1)} \,\,\,$A:D(A)\subset H\to H$ and $B:D(B)\subset H\to H$ are unbounded, monotone,\linebreak 
\hspace*{14mm}self-adjoint, linear operators with compact resolvents.  
     
\vspace{2mm}\noindent
Therefore, there are sequences 
$\{\lambda_j\}$, $\{\lambda'_j\}\,$ and $\,\{e_j\}$, $\{e'_j\}$ 
of eigenvalues and of corresponding eigenfunctions such that 
\begin{align}
\label{eigen}
& A e_j = \lambda_j e_j, \quad
B e'_j = \lambda'_j e'_j, \quad \hbox{with} \quad
 (e_i,e_j) = (e'_i,e'_j) = \delta_{ij}
\quad \forall\,i,j\in\enne,
\\[1mm]
\label{eigenva}
& 0 \leq \lambda_1 \leq \lambda_2 \leq \dots , \quad
  0 \leq \lambda'_1 \leq \lambda'_2 \leq \dots, \quad \hbox{with} \quad
  \lim_{j\to\infty} \lambda_j
  = \lim_{j\to\infty} \lambda'_j
  = + \infty,
\\[1mm]
\label{eigenfu}
& \hbox{$\{e_j\}$\, and \,$\{e'_j\}$\, are complete systems in $H$}.
  \end{align}
As a consequence, we can define the powers of \,$A$\, and \,$B$\, for arbitrary 
positive real exponents: we have, for $\rho>0$,
\Bsist
  && \VA \rho := D(A^\rho)
  = \Bigl\{ v\in H:\ \somma j1\infty |\lambda_j^\rho (v,e_j)|^2 < +\infty \Bigr\}
  \aand
  \label{domAr}
  \\[-3mm]
  && A^\rho v = \somma j1\infty \lambda_j^\rho (v,e_j) e_j
  \quad \hbox{for $v\in\VA \rho$},
  \label{defAr}
\Esist
the series being convergent in the strong topology of~$H$.
By endowing $\VA \rho$ with the graph norm, i.e., setting
\Beq
  (v,w)_{\VA \rho} := (v,w) + (A^\rho v,A^\rho w)
  \aand
  \norma v_{\VA \rho} := (v,v)_{\VA \rho}^{1/2}
  \quad \hbox{for $v,w\in\VA \rho$},
  \label{normAr}
\Eeq
we obtain a Hilbert space. 
In the same way, 
we define the power $B^\sigma$ for every $\sigma>0$,
starting from \eqref{eigen}--\eqref{eigenfu} for~$B$.
We therefore set $\,\VB\sigma := D(B^\sigma),$ 
with the norm $\norma\cpto_{\VB \sigma}$
associated with the inner product
\begin{equation}
\label{normBs}
(v,w)_{\VB \sigma} := (v,w) + (B^\sigma v,B^\sigma w)
\quad \hbox{for $v,w\in \VB\sigma$.}
\end{equation}
Since $\lambda_j\geq0$ \juerg{and $\lambda_j'\ge 0$} for every~$j$, one immediately deduces from the definitions that
$\,A^\rho:\VA\rho\subset H\to H\,$ and $\,B^\sigma:\VB \sigma\subset H\to H\,$ 
 are maximal monotone operators. 
Moreover, it is clear that, for every $\rho_1,\,\rho_2>0$, we have the Green type formula
\Beq
  (A^{\rho_1+\rho_2} v,w)
  = (A^{\rho_1} v, A^{\rho_2} w)
  \quad \hbox{for every $v\in\VA{\rho_1+\rho_2}$ and $w\in\VA{\rho_2}$},
  \label{Green}
\Eeq
and that a similar relation holds  for~$\,B$.
Due to these properties, we can define proper extensions of the operators that
allow values in dual spaces.
In particular, we can write variational formulations of \eqref{sstheta} and \eqref{ssphi}.
It is convenient to use the notations
\Beq
  \VA{-\rho} := (\VA \rho)^* , \quad
  \VB{-\sigma} := (\VB\sigma)^*,
    \quad \hbox{for $\rho>0$ and $\sigma>0$}.
  \label{negspace}
\Eeq
Then, we have that
\Beq
\label{ext1}
  A^{2\rho} \in \calL(V_A^\rho, V_A^{-\rho}) , \quad
  B^{2\sigma} \in \calL(V_B^\sigma, V_B^{-\sigma}),
\Eeq
as well as
\Beq
\label{ext2}
  A^\rho \in \calL(H, V_A^{-\rho}) , \quad
  B^\sigma \in \calL(H, V_B^{-\sigma}).   
\Eeq
Here, we identify $H$ with a subspace of $V_A^{-\rho}$
in the usual way, i.e., such that
\Beq
  \< v,w >_{V_A^{\rho}} = (v,w)
  \quad \hbox{for every $v\in H$ and $w\in V_A^\rho$}.
  \label{identify}
\Eeq
Analogously, we have $H\subset\VB{-\sigma}$ and use corresponding notations.
Observe also that the following embeddings are continuous and compact:
\begin{align}
\label{embed1}
&\VA{\rho_1+\rho_2}\subset \VA{\rho_1}\subset H,\quad \VB{\sigma_1+\sigma_2}\subset 
\VB{\sigma_1}\subset H,\quad\mbox{for $\rho_1>0$, $\rho_2>0$\, and \,$\sigma_1>0$, $\sigma_2>0$,}
\\[0.5mm] 
\label{embed2}
&\VA \rho \subset H \subset \VA{-\rho}, \quad \VB \sigma\subset H\subset \VB{-\sigma},
\quad\mbox{for $\,\rho>0\,$ and \,$\sigma>0$.}
\end{align}

From now on, we generally assume for the nonlinear functions entering \eqref{sstheta} and \eqref{ssphi}:

\vspace{2mm}\noindent
{\bf (F1)} \,\,\,$F=F_1+F_2$, where $\,F_1:\erre\to [0,+\infty]\,$ is convex and lower semicontinuous
with \linebreak \hspace*{12mm}$\,F_1(0)=0$. Moreover, there are constants $c_1>0$, $c_2>0$, such that 
\Beq
\label{coerc}
F(s)\ge c_1s^2-c_2 \quad\forall\,s\in\erre.
\Eeq

\vspace{2mm}\noindent
{\bf (F2)} \,\,\,There are $\,r_-,r_+\,$ with $\,-\infty\le r_-<0<r_+\le +\infty$\, such that
$\,F_1\in C^3(r_-,r_+)$, and \linebreak \hspace*{13mm}it holds $\,F_1'(0)=0$.

\vspace{2mm}\noindent
{\bf (F3)} \,\,\,$F_2\in C^3(\erre)$, and $F_2'$ is Lipschitz continuous on $\erre$ with Lipschitz
constant $L>0$.

\vspace{2mm}\noindent
{\bf (F4)} \,\,\,$\ell\in C^2(\erre)$, and $\ell^{(i)}\in L^\infty(\erre)$ for $0\le i\le 2$.

\vspace{1mm}
\Brem 
It is worth noting that all of the potentials \eqref{regpot}--\eqref{obspot}
satisfy the general conditions {\bf (F1)}--{\bf \pier{(F3)}}, 
where $\,D(F_1)=D(\partial F_1)=\erre$ for $F=F_{\rm reg}$, 
while 
\gianni{$D(F_1)=[-1,1]$ and $D(\partial F_1)=(-1,1)$ for $F=F_{\rm log}$,
and $D(F_1)=D(\partial F_1)=[-1,1]$ for $F=F_{\rm 2obs}$}. 
Here, and throughout this
paper, we denote 
by $\,D(F_1)\,$ and $\,D(\partial F_1)\,$ the effective domains
of $\,F_1\,$ and of its subdifferential $\,\partial F_1$, respectively.
We notice that $\partial F_1$ is a maximal monotone graph in $\erre\times\erre$
and use the same symbol $\partial F_1$ for the maximal monotone operators induced in $L^2$ spaces.
Moreover, for $r\in D(\partial F_1)$, we use the symbol $\partial F_1^{\rm o}(r)$\, for the element 
of $\,\partial F_1(r)\,$ having minimal modulus. 
If, however, $\partial F_1$ is single-valued
(which is the case if $(r_-,r_+)=\erre$), 
then we denote the sole element of the singleton
$\partial F_1(r)$ by~$\,F_1'(r)$. 
We also remark that {\bf (F3)} implies that $F_2'$ grows at most linearly on $\erre$, while 
$F_2$ grows at most quadratically.  
\Erem

\vspace{3mm}
For the other data of the state system, we postulate:  

\vspace{2mm}\noindent
{\bf (A2)} \,\,\,$\rho$ \,and \, $\sigma$\, are fixed positive real numbers.

\vspace{2mm}\noindent
{\bf (A3)}  \,\,\,$\theta_0\in \VA \rho$, $\,\phi_0\in V_B^{2\sigma}$, and there are constants
$\,r_{0-}, r_{0+}\,$ such that
\begin{equation}
\label{sepini}
r_-<r_{0-}\,\le\,\phi_0\,\le r_{0+} <r_+ \quad\mbox{a.e.\ in }\,\Omega. 
\end{equation}

\vspace{2mm}\noindent
{\bf (A4)}  \,\,\,The embeddings $\,V_A^\rho\subset L^4(\Omega)\,$ and $\,V_B^\sigma\subset \Lx 4\,$
are continuous.

\vspace{1mm}
\Brem
\label{RemA4}
If, for instance, \gianni{$A=-\Delta$} with domain $\Hdue\cap\Hunoz$ 
(thus, with homogeneous Dirichlet conditions, but similarly for zero boundary conditions of
Neumann or third kind),
then $\VA\rho \subset H^{2\rho}(\Omega)$; 
it then follows from \eqref{sobslo3} that \pier{(the first embedding in)}~{\bf (A4)} 
holds true if $\rho\ge 3/8$. 
Likewise, we have in this case $\,\VA\rho \subset \Lx 6\,$ for
\gianni{$\,\rho\ge 1/2$} as well as $\VA\rho\subset C^0(\overline\Omega)$ provided that $\,\rho>3/4$.
\Erem

\vspace{2mm} For the data entering the cost functional \eqref{cost} and the admissible set
$\uad$ defined in \eqref{uad} we generally assume:

\vspace{2mm}\noindent
{\bf (A5)} \,\,\,$\theta_\Omega,\phi_\Omega\in\Ldue$, \,\,$\theta_Q,\phi_Q\in L^2(Q)$, 
\,$\,u_{min},u_{max}\in L^\infty(Q)\,$ satisfy $\,u_{min}\le u_{max}$\linebreak
\hspace*{13mm} a.e.\ in $\,Q$.

\vspace{2mm} Finally, \pier{once and for all we} fix some open  and bounded ball in $L^\infty(Q)$ that contains
the admissible set.

\vspace{2mm}\noindent
{\bf (A6)} \,\,\,$R>0$ is a constant such that $\,\uad\subset {\cal U}_R:=\left\{u\in L^\infty(Q):
\,\|u\|_{L^\infty(Q)}<R\right\}$. 

\vspace{3mm}
At this point, we are in a position to make use of \eqref{Green} and its analogue for $B$ 
to give a weak formulation of the state system \eqref{sstheta}--\eqref{ssini} and to introduce our notion of solution.  
In particular, we present \eqref{ssphi} in the form of a variational inequality.
We look for a pair $(\theta,\phi)$ satisfying
\begin{align}
\label{regtheta}
&\theta\in H^1(0,T;H)\cap L^\infty(0,T;\VA\rho)\cap L^2(0,T;\VA{2\rho}),\\[0.5mm]
\label{regphi}
&\phi\in W^{1,\infty}(0,T;H)\cap H^1(0,T;\VB\sigma),\\[0.5mm]
\label{L1Q}
&F_1(\phi)\in L^1(Q), 
\end{align}
and solving the system
\begin{align}
\label{ssvar1}
&\dt\theta+\ell(\phi)\dt\phi+A^{2\rho}\theta=u\quad\mbox{a.e.\ in }\,Q,\\
\label{ssvar2}
&(\dt\phi(t),\phi(t)-v)+(B^\sigma\phi(t),B^\sigma(\phi(t)-v))+\iO F_1(\phi(t))
\,+\,(F_2'(\phi(t)),\phi(t)-v)\nonumber\\          
&\le\,(\pier{\ell(\phi(t))}\theta(t),\phi(t)-v)\,+\iO F_1(v)\quad\mbox{for a.e.\ }\,t\in (0,T) \,\mbox{ and 
every }\,v\in\VB\sigma,\\
\label{ssvar3}
&\theta(0)=\theta_0,\quad \phi(0)=\phi_0.
\end{align}
Here, it is understood that $\,\,\iO F_1(v)=+\infty\,$ whenever $\,F_1(v)\not\in L^1(\gianni\Omega)$. 
We follow a similar rule for expressions of the type $\,\iint_QF_1(v)\,\,$ whenever
$\,v\in L^2(Q)\,$ but $F_1(v)\not\in L^1(Q)$.
 
We notice that \eqref{ssvar2} is equivalent to its time-integrated variant, that is,
\begin{align}
  &  \ioT \bigl( \dt\phi(t) , \phi(t) - v(t) \bigr) \, dt 
  + \ioT \bigl( B^\sigma\phi(t) , B^\sigma(\phi(t)-v(t)) \bigr)\,dt 
  \nonumber
  \\ 
  & \quad {}
  + \juerg{\iint_Q} F_1(\phi) 
  + \ioT \bigl( F_2'(\phi(t)) , \phi(t)-v(t) \bigr) \, dt 
  \nonumber
  \\
  & \leq \ioT \bigl( \ell(\phi(t))\theta(t) , \phi(t)-v(t) \bigr)\,dt
  + \juerg{\iint_Q} F_1(v)
  \quad \hbox{for all $v\in\L2{\VB\sigma}$}.
  \label{intssvar2}
\end{align}
Similarly, \eqref{ssvar1} is equivalent to a corresponding time-integrated version with test functions $v\in
L^2(0,T;\VA\rho)$.

We have the following well-posedness result (cf.\ \cite[Thm.~2.10]{CG}).

\Bthm
\label{FromCG}
Let the assumptions {\bf (F1)}--{\bf (F4)}, {\bf (A1)}--{\bf (A4)}, and {\bf (A6)}
be fulfilled. 
Then the problem \eqref{ssvar1}--\eqref{ssvar3} has for every $u\in{\cal U}_R$
a unique  solution $(\theta,\phi)$ satisfying \eqref{regtheta}--\eqref{L1Q}. 
Moreover,
there is a constant $K_1>0$, which depends only on $R$ and the data of the state system,
such that 
\begin{align}
\label{ssbound1}
&\|\theta\|_{H^1(0,T;H)\cap L^\infty(0,T;\VA\rho)\cap L^2(0,T;\VA {2\rho})}
 \,+\, \|\phi\|_{W^{1,\infty}(0,T;H)\cap H^1(0,T;\VB \sigma)} 
\,+\,\juerg{\iint_Q} F_1(\phi)\,\le\,K_1,
\end{align}  
whenever $(\theta,\phi)$ solves \eqref{ssvar1}--\eqref{ssvar3} for some $u\in{\cal U}_R$.
\Ethm

\Bdim
Owing to the assumptions {\bf (F2)} and \eqref{sepini}, we have that $\,\partial F_1^{\rm o}(\phi_0)
=F_1'(\phi_0)\in H$, and thus all of the conditions for the application of \cite[Thm.~2.10]{CG} are fulfilled.
\Edim

\vspace{2mm}
By virtue of Theorem~\ref{FromCG}, the {\em control-to-state operator} 
\begin{equation}\label{defS}
{\cal S}: \,{\cal U}_R\ni u\mapsto {\cal S}(u):=(\theta,\phi) 
\end{equation}
is well defined and bounded as a mapping from ${\cal U}_R\subset L^\infty(Q)$ into the Banach space specified
by the regularity properties \eqref{regtheta}, \eqref{regphi}.
In the following, we look for conditions that guarantee that $\,{\cal S}\,$ is Fr\'echet differentiable between 
suitable Banach spaces. 
To this end, we make use of \gianni{the following global} boundedness assumption which has proved
to be very useful in the framework of Cahn--Hilliard type systems with fractional operators 
(cf.\ the recent works \cite{CGS19, CGS21, CGS21bis,CGS25}):

\vspace{3mm}\noindent
{\bf (GB)} \quad There are constants $a_R,b_R$ 
such that
\begin{equation}
\label{(GB)}
r_-<a_R\le\varphi\le b_R< r_+
\quad\mbox{a.e.\ in }\, Q
\end{equation}
\hspace*{17mm}whenever $(\vartheta,\varphi)={\cal S}(u)$ for some $u\in{\cal U}_R$.

\vspace{2mm}
The condition {\bf (GB)} is rather restrictive and has to be verified from case to case. 
As a rule, it cannot be satisfied for potentials of indicator function type like
$F_{\rm 2obs}$. It can, however, be satisfied for regular potentials like $F_{\rm reg}$ 
(see the case (ii) below) 
and singular potentials like $F_{\rm log}$ 
(see the case (i)  below). 
Indeed, we have the following result, 
\gianni{whose assumptions are commented in the next Remark~\ref{Rempartin}.}

\Blem
\label{SuffGB}
\pier{Assume that} the conditions {\bf (A1)}--{\bf (A4)} and {\bf (F1)}--{\bf (F4)} are satisfied \pier{and, in addition, that}
\begin{align}
\label{partin}
&\psi(v)\in H\,\, \mbox{ and }\,\,(B^{2\sigma}v,\psi(v))\ge 0 \quad\mbox{for every $\,v\in V_B^{2\sigma}$\, and every
monotone}\nonumber\\
&\quad \mbox{and Lipschitz continuous mapping $\,\psi:\erre\to\erre\,$ vanishing at the origin}.
\end{align}
In addition, assume that
\begin{equation}\label{asymp}
\lim_{r\to r_-} F_1'(r)=-\infty,\quad\lim_{r\to r_+} F_1'(r)=+\infty.
\end{equation}
Then {\bf (GB)} is satisfied in any of the following situations:\\[1mm]
\hbox to 2em{{\rm (i)}\hfil} $A=-\Delta$\, with zero Neumann or Dirichlet boundary conditions,
  and \,$\rho>\frac34$ or $\,\rho=\frac 12$.
\\[1mm]
\hbox to 2em{{\rm (ii)}\hfil} $(r_-,r_+)=\erre$, \gianni{$B^{2\sigma}=B=-\Delta$ with zero Dirichlet or Neumann boundary conditions 
\hfil\break\hbox to 2em{\hfil} and $\phi_0\in D(B)$}.
\Elem

\Bdim 
Suppose that $(\theta,\phi)={\cal S}(u)$ for some $u\in{\cal U}_R$. 
Assume first that the assumptions of (i) are satisfied.  
In the following, we denote by $C_i>0$, $i\in\enne$, constants that depend
only on $R$ and the data. We have, owing to  \eqref{sobslo2}, that
\pier{(cf.~also Remark~\ref{RemA4})}
$V_A^\rho\subset H^{2\rho}(\Omega)\subset L^\infty(\Omega)$ if $\rho>3/4$. 
Hence, by \eqref{ssbound1} \pier{it turns out that}
\begin{equation}\label{yeah}
\|\theta\|_{L^\infty(Q)}\,\le\,C_1
\end{equation}
in this case. 
On the other hand, if $\rho=1/2$, then $\theta$
solves a standard linear parabolic problem with right-hand side $\,u-\ell(\phi)\dt\phi$, which is bounded in 
$L^\infty(0,T;H)$ for $u\in{\cal U}_R$.
Then the validity of \eqref{yeah} follows from standard results on linear parabolic problems 
(see, e.g., \cite[Chap.~7]{LSU}). 
Hence, in both cases, \pier{we have that}
$\|\ell(\phi)\theta\|_{L^\infty(Q)}\,\le\,C_2$ \gianni{since $\ell$ is bounded}, 
and the validity of {\bf(GB)}
with $a_R, b_R$ satisfying $\,r_-<a_R\le r_{0-}\le r_{0+}\le b_R<r_+$ follows from \gianni{the assumptions \accorpa{partin}{asymp}} 
as in the proof of \cite[Thm.~\gianni{2.4}]{CGS25}.
 
Now\pier{, let} the assumptions of (ii) \pier{be} fulfilled.
Then\pier{, we remark that $(r_-,r_+)=\erre$ excludes} singular potentials like $F_{\rm log}$ \gianni{and we have to prove that $\phi$ is bounded in $\LQ\infty$ uniformly with respect to $u\in\calU_R$.
Let}, 
for $\lambda>0$, $F_{1,\lambda}'$ denote the Moreau--Yosida
approximation of $F_1'$ at the level $\lambda$. 
It is well known (see, e.g., \cite{Brezis}) 
that in this special case, where the subdifferentials are \pier{single-valued} \gianni{and $F_1'(0)=0$}, the following
conditions are satisfied:
\begin{align}\label{Yosi}
&F'_{1,\lambda} \mbox{\, is globally Lipschitz continuous on $\,\erre$, \,\,\gianni{$F'_{1,\lambda}(0)=0$},\, and it holds}\nonumber\\
&|F'_{1,\lambda}(r)|\le |F_1'(r)| \,\,\mbox{ and }\,\,\lim_{\lambda\searrow0}\,F'_{1,\lambda}(r)=F_1'(r)
\quad\mbox{for all }\,r\in\erre.  
\end{align} 
In the proofs of \gianni{\cite[Prop.~2.4 and Prop.~2.9]{CG}}
(for details, see \cite[Sect.~5]{CG}) 
it has been shown, \gianni{using \eqref{asymp} and a special case of~\eqref{partin},} 
that there is some $\Lambda>0$ such that for every $\lambda\in (0,\Lambda]$ the \gianni{general} system
\gianni{%
\begin{align}
&\dt\thetal + \ell(\phil)\dt\phil+A^{2\rho}\thetal\,=\,u\quad\mbox{a.e.\ in }\,Q,
\label{MY1}
\\[1mm]
&\dt\phil+B^{2\sigma}\phil+\F1l'(\phil)+F_2'(\phil)\,=\,\ell(\phil)\thetal\quad\mbox{a.e.\ in }\, Q,
\label{MY2}
\\[1mm]
&\thetal(0)=\theta_0, \quad \phil(0)=\phi_0 \quad\mbox{a.e.\ in }\,\Omega,
\label{MY4}
\end{align}
}%
has for every $u\in {\cal U}_R$ a unique solution pair $(\thetal,\phil)$ such that
\begin{align}
\label{boundl}
\|\thetal\|_{H^1(0,T;H)\cap L^\infty(0,T;V_A^\rho)\cap L^2(0,T;V_A^{2\rho})}
\,+\,\|\phil\|_{W^{1,\infty}(0,T;H)\cap H^1(0,T;\VB\sigma)\cap L^2(0,T;\VB {2\sigma})}\,\le\,M_1,
\end{align}
where, here and in the following, $M_i>0$, $i\in\enne$, denote constants that may depend on $R$ and 
on the data of the system,
but not on $\lambda\in (0,\Lambda]$.
\gianni{It \juerg{was} then shown in \cite{CG} that $(\thetal,\phil)$ converges to $(\theta,\phi)$ in a suitable topology.}
\gianni{We repeat here a part of the argument and use \eqref{partin} with $v=\phil(t)$ and $\psi=F'_{1,\lambda}$.}
We test \eqref{MY2} by $F'_{1,\lambda}(\phil(t))$ to obtain for almost every $t\in(0,T)$ the \juerg{identity}
\begin{align}
\label{idl}
& \gianni{\bigl( B^{2\sigma}\phil(t),F'_{1,\lambda}(\phil(t))\bigr)
  + \iO |F'_{1,\lambda}(\phil(t))|^2}
\nonumber\\
&=\iO F_{1,\lambda}'(\phil(t))\bigl(\ell(\phil(t))\thetal(t)-F_2'(\phil(t))-\dt\phil(t)\bigr),
\end{align}
where the \gianni{first} summand on the left-hand side is nonnegative and, by virtue of \eqref{boundl} and the 
general assumptions for the nonlinearities, the right-hand side is bounded by an expression of the form
$$
\frac 12 \,\iO\left|F'_{1,\lambda}(\phil(t))\right|^2 \,+\,M_2,
$$
whence we obtain that
$$
\bigl\|F'_{1,\lambda}(\phil)\bigr\|_{L^\infty(0,T;H)}\,\le\,M_3.
$$
\gianni{Thanks to our assumption on~$B^{2\sigma}$ it follows} that $\,\phil\,$ solves a
standard linear parabolic initial-boundary value problem, where
the right-hand side $\,\ell(\phil)\thetal-F'_{1,\lambda}(\phil) -F'_2(\phil)\,$ is bounded in $L^\infty(0,T;H)$,
independently of $\lambda\in (0,\Lambda]$. 
\gianni{Moreover, $\phi_0\in D(B)\subset\Linfty$.}
Applying the classical results of \cite[Chap.~7]{LSU}, we therefore can infer
that $\,\|\phil\|_{L^\infty(Q)}\,$ is bounded independently of $\lambda\in (0,\Lambda]$. 
Hence, $\phil\to\phi$\, weakly-star
in $L^\infty(Q)$, and the lower semicontinuity of norms yields the assertion.
\Edim

\Brem
\label{Rempartin}
\gianni{The assumptions on $B^{2\sigma}$ made in (ii), 
and the condition \eqref{partin} we also used in the second part of the proof, are not in contradiction with each other.
In fact, the former implies the latter.
Indeed, in this case,
$\VB\sigma=\Hunoz$ or $V_B^\sigma = H^1(\Omega)$, and therefore it holds for every monotone and 
Lipschitz continuous function $\,\psi\,$ vanishing at
the origin that for every $v\in V_B^{2\sigma}\subset\Hdue$ we have $\psi(v)\in\Huno$, as well~as
$$
(B^{2\sigma}v,\psi(v))\gianni{{}=(-\Delta v,\psi(v))}=\iO\psi'(v)\,|\nabla v|^2 \,\ge 0.
$$ 
Moreover, in both (i) and~(ii), the Laplacian
can be replaced by more general second-order elliptic operators in divergence form 
with smooth coefficients complemented with more general zero boundary conditions.
Furthermore, regarding $A$ in~(i), even higher order operators can be considered 
provided that the assumptions on $\rho$ are modified accordingly.
For instance, one can take the plate operator $A=\Delta^2$,  assuming as domain 
the set of $v\in\Hx4$ satisfying suitable boundary conditions.
Two possibilities are $v=\dn v=0$ and $\dn v=\dn\Delta v=0$.
In both cases, $\VA\rho\subset\Hx{4\rho}$, and then the condition $\VA\rho\subset\Linfty$ 
used in the above proof is satisfied if $\rho>3/8$.} 
\Erem

\Brem
\label{FromGB}
If condition {\bf (GB)} is fulfilled, then it follows from the assumptions {\bf (F2)}--{\bf (F4)}
that, by possibly taking a larger constant $K_1>0$, 
it holds the global bound 
\Beq
\max_{0\le i\le 3}\,\|F_1^{(i)}(\phi)\|_{L^\infty(Q)}\,+\,\max_{0\le i\le 3}\,\|F_2^{(i)}(\phi)\|_{L^\infty(Q)}
  \,+\,\max_{0\le i\le 2}\,\|\ell^{(i)}(\phi)\|_{L^\infty(Q)}\,\le\,K_1
  \label{ssbound2}
\Eeq
whenever $(\theta,\phi)={\cal S}(u)$ for some $u\in{\cal U}_R$.
\Erem

\vspace{3mm} 
The next step in our analysis is to show that under the condition {\bf (GB)}
a rather strong stability estimate holds true for the solutions to the state system,
which constitutes an important preparation for the later proof of Fr\'echet differentiability.
We have, however, to make a further assumption:

\vspace{2mm}\noindent
{\bf (A7)} \,\,\,$\VB\sigma\cap L^\infty(\Omega)$ is dense in $\VB\sigma$.

The condition {\bf (A7)} is, for example, fulfilled if $\VB\sigma$ coincides with one of the
Sobolev--Slobodeckij spaces $H^s(\Omega)$ for $s>0$. 
\gianni{We combine it with {\bf(GB)} to prove the lemma below. 
Similar results \juerg{were established in} \cite[Prop.~2.4 and Prop.~2.9]{CG} under an assumption close to~\eqref{partin}.}

\gianni{%
\Blem
\label{FromA7}
Suppose that the conditions {\bf (F1)}--{\bf (F4)}, {\bf (A1)}--{\bf (A4)}, {\bf (A6)}--{\bf (A7)}, and
{\bf (GB)} are satisfied.
Then, \juerg{for every $u\in\calU_R$,} the solution $(\theta,\phi)$ to the state system satisfies the variational equation
\begin{align}
  \label{ssvar2n}
  &(\dt\phi(t),v)+(B^\sigma\phi(t),B^\sigma v)+(F'(\phi(t)),v)\,=\,(\ell(\phi(t))\theta(t),v)\nonumber\\
  &\quad\mbox{for a.e.\ $t\in\gianni{(0,T)}$ and all }\,v\in\VB\sigma.
\end{align}
Moreover, by possibly enlarging the constant $K_1$ that appears in \eqref{ssbound1} and~\eqref{ssbound2},
we have the estimate
\Beq
  \norma\phi_{\W{1,\infty}H\cap\H1{\VB\sigma}\cap\L\infty{\VB{2\sigma}}} \leq K_1 \,.
    \label{moreregphi}
\Eeq
In particular, $(\theta,\phi)$ is a strong solution to the state system.
\Elem
}%

\Bdim
By recalling {\bf(GB)}, we set $\delta:=\min\{a_R-r_-,r_+-b_R\}$.
Take now an arbitrary $w\in\VB\sigma\cap\Linfty$,
and let $\eps_0>0$ be such that $\eps_0\norma w_{\Linfty}\leq \delta/2$.
Then, $v:=\phi(t)+\eps w$ is for every $\eps\in(0,\eps_0)$ 
an admissible test function in~\eqref{ssvar2}, and $F_1(v)\in\Luno$.
By using it and then dividing by~$-\eps$, we obtain (\aat)
\Bsist
  && \bigl( \dt\phi(t),w \bigr)
  + \bigl( B^\sigma \phi(t) , \juerg{B^\sigma}w \bigr)
  + \iO \frac {F_1(\phi(t)+\eps w) - F_1(\phi(t))} \eps
  \non
  \\
  && \quad {}
  + \bigl( F_2'\pier{(\phi(t))}, w \bigr)
  \geq \bigl( \ell(\phi(t)) \theta(t) , w \bigr).
  \non
\Esist
Since $r_-+\delta/2\leq\phi(t)+\eps w\leq r_+-\delta/2$ for every $\eps\in(0,\eps_0)$,
and since $\phi(t)+\eps w$ converges to $\phi(t)$ in the strong topology of $\VB\sigma\cap\Linfty$ as $\eps\searrow0$,
we immediately deduce that
\Beq
  \bigl( \dt\phi(t),w \bigr)
  + \bigl( B^\sigma \phi(t) , w \bigr)
  + \bigl( F_1'\juerg{(\phi(t))} , w \bigr)
  + \bigl( F_2'\juerg{(\phi(t))} , w \bigr)
  \geq \bigl( \ell(\phi(t)) \theta(t) , w \bigr).
  \non
\Eeq
By changing $w$ into $-w$, we obtain the opposite inequality, and thus equality.
Finally, by accounting for~{\bf(A7)}, we can remove the boundedness assumption on the test function.

Let us come to~\eqref{moreregphi}.
A~part of it is already given by~\eqref{ssbound1}.
To derive the maximal space regularity,
we notice that \eqref{ssvar2n} can be written~as
\Beq
  B^{2\sigma} \phi = g := \ell(\phi) \theta - \dt\phi - F'(\phi)
  \quad \hbox{in $\VB{-\sigma}$, \aet},
  \non
\Eeq
and that $g$ is bounded in $\L\infty H$ by a constant that only depends on $R$ and the data
thanks to \eqref{ssbound1} and~\eqref{ssbound2}.
\Edim
 
We \gianni{derive} the following result.

\Bthm
\label{Contdep}
Suppose that the conditions {\bf (F1)}--{\bf (F4)}, {\bf (A1)}--{\bf (A4)}, {\bf (A6)}--{\bf (A7)}, and
{\bf (GB)} are satisfied. 
Then there is a constant $K_2>0$, which depends only on $R$ and
the data of the state system, such that the following holds true: 
whenever $(\theta_i,\phi_i)
={\cal S}(u_i)$, $i=1,2$, for some controls \,$u_1,u_2\in{\cal U}_R$, then, for every $t\in (0,T]$, 
\begin{align}\label{ssbound3}
&\|\theta_1-\theta_2\|_{H^1(0,t;H)\cap L^\infty(0,t;\VA\rho)}\,+\,
\|\phi_1-\phi_2\|_{W^{1,\infty}(0,t;H)\cap H^1(0,t;\VB\sigma)}\nonumber\\
&\le\,K_2\,\|u_1-u_2\|_{L^2(0,t;H)}.
\end{align}
\Ethm  
    
\Bdim 
Since {\bf (GB)} is fulfilled, the global bounds \eqref{ssbound1} and \eqref{ssbound2} are 
satisfied for $(\theta_i,\phi_i)$, $i=1,2$.
\gianni{Moreover, by Lemma~\ref{FromA7}, we can replace 
the variational inequality \eqref{ssvar2} by the variational equation~\eqref{ssvar2n}.}
Now let $\theta:=\theta_1-\theta_2$,
$\phi:=\phi_1-\phi_2$, and $\,u:=u_1-u_2$. 
Then it is easily seen that $(\theta,\phi)$ is a strong solution to the system
\begin{align}
\label{diff1}
&\dt\theta+A^{2\rho}\theta+(\ell(\phi_{1})-\ell(\phi_{2}))\dt\phi_{1}+\ell(\phi_{2})\dt\phi\,=\,u
\quad\mbox{in }\,Q,\\[0.5mm]
\label{diff2}
&\dt\phi+B^{2\sigma}\phi+F'(\phi_{1})-F'(\phi_{2})\,=\,(\ell(\phi_{1})-\ell(\phi_{2}))\theta_{1}+\ell(\phi_{2})\theta
\quad\mbox{in }\,Q,\\[0.5mm]
\label{diff3}
&\theta(0)=0,\quad \phi(0)=0,\quad\mbox{in }\,\Omega.
\end{align} 

To begin with, we test \eqref{diff1} by $\theta$ and \eqref{diff2} by $\,\dt\phi\in L^2(0,T;\VB\sigma)$, add the resulting 
equations and integrate over $\Omega\times(0,t)$, where $t\in(0,T)$. 
Adding the same term $\gianni{\frac12\,\norma{\phi(t)}^2={}}\int_0^t\!\iO\phi\,\dt\phi$ 
to both sides of the resulting identity \gianni{and noting a cancellation}, we arrive at the equation
\begin{align}
\label{diff5}
&\frac 12 \,\|\theta(t)\|^2\,+\,\frac 12\,\|\phi(t)\|_{\VB \sigma}^2\,+\int_0^t\!\!\iO|A^\rho\theta|^2
+\int_0^t\!\!\iO|\dt\phi|^2\nonumber\\
&=\,\int_0^t\!\!\iO u\theta\,-\int_0^t\!\!\iO\theta(\ell(\phi_1)-\ell(\phi_2))\dt\phi_1\,+\int_0^t\!\!\iO
(\ell(\phi_1)-\ell(\phi_2))\,\theta_1\,\dt\phi\nonumber\\
&\hspace*{5mm}- \int_0^t\!\!\iO\bigl(F'(\phi_1)-F'(\phi_2)\bigr)\dt\phi\,
+\int_0^t\!\!\iO\phi\,\dt\phi\,=:\,\sum_{j=1}^5 I_j\,,
\end{align}
with obvious notation. We estimate the terms on the \rhs\ individually, using the Young and H\"older inequalities,
the embedding conditions of {\bf (A4)}, 
as well as the global bounds \eqref{ssbound1} and \eqref{ssbound2}, repeatedly without further reference.
In this process, $C_i$, $i\in\enne$, denote constants that depend only on $R$ and the data of the state system. 
Clearly, we have
\begin{equation}\label{diff6}
|I_1|\,\le\,\frac 12\int_0^t\!\!\iO\theta^2\,+\,\frac 12\int_0^t\iO|u|^2\,.
\end{equation}
Moreover, for every $\delta>0$ (which is yet to be specified) it follows that                                             
\begin{align}\label{diff7}
|I_2|\,&\le\,C_1\int_0^t\|\theta(s)\|_{\Lx 4}\,\|\phi(s)\|_{\Lx 4}\,\|\dt\phi_1(s)\|\,ds\nonumber\\
&\le\,\delta\int_0^t\|\theta(s)\|_{\VA \rho}^2\,ds\,+\,\frac{C_2}\delta\int_0^t\|\dt\phi_1(s)\|_{V_B^\sigma}^2\,
\|\phi(s)\|_{V_B^\sigma}^2\,ds\,.
\end{align}
In  addition, \pier{we see that}
\begin{align}\label{diff8}
|I_3|\,&\le\,C_3\int_0^t\|\theta_1(s)\|_{\Lx 4}\,\|\phi(s)\|_{\Lx 4}\,\|\dt\phi(s)\|\,ds\nonumber\\
&\le\,\delta\int_0^t\!\!|\dt\phi|^2\,+\,\frac{C_4}\delta\int_0^t\|\theta_1(s)\|^2_{\VA\rho}\,
\|\phi(s)\|^2_{\VB\sigma}\,ds\,.
\end{align}
Finally, owing to \pier{the Lipschitz continuity of $F'$ in $[a_R, b_R]$, 
it turns out that}
\begin{equation}\label{diff9}
|I_4|+|I_5|\,\le\,\delta\int_0^t\!\!\iO|\dt\phi|^2\,+\,\gianni{\frac{C_5}\delta}\int_0^t\!\!\iO|\phi|^2\,.
\end{equation}
Now observe that the mapping $\,s\mapsto\|\theta_1(s)\|_{\VA\rho}^2 +\|\dt\phi_1(s)\|^2_{\VB\sigma}\,\,$
belongs by \eqref{ssbound1} to $L^1(0,T)$. 
Hence, choosing $\delta=1/4$, we can infer from
Gronwall's lemma that for every $t\in (0,T)$ we have the estimate
\begin{equation}\label{diff10}
\|\theta\|_{L^\infty(0,t;H)\cap L^2(0,t;\VA\rho)}\,+\,\|\phi\|_{H^1(0,t;H)\cap L^\infty(0,t;\VB\sigma)}
\,\le\,C_6\,\|u\|_{L^2(0,t;H)}\,.
\end{equation}

In the next estimate, we argue formally, noting that the arguments can be made rigorous by using, e.g.,
finite differences in time and the fact that $\theta(0)=\phi(0)=0$. 
Indeed, we formally differentiate \gianni{\eqref{diff2}} with respect to time to obtain the identity
\begin{align}\label{diff11}
&\partial_{tt}\phi+B^{2\sigma}\dt\phi+(F''(\phi_1)-F''(\phi_2))\dt\phi_1+F''(\phi_2)\dt\phi\nonumber\\
&=\,(\ell'(\phi_1)-\ell'(\phi_2))\theta_1\,\dt\phi_1+\ell'(\phi_2)\theta_1\,\dt\phi
+\dt\theta_1(\ell(\phi_1)-\ell(\phi_2))\nonumber\\
&\quad\,\pier{{}+\theta\,\ell'(\phi_2)\dt\phi_2 + \ell(\phi_2)\dt\theta}\,.
\end{align}
\pier{Let us add $\dt\phi$ to both sides of \eqref{diff11}. Then, we} (formally) test \gianni{\eqref{diff1}} by $\dt\theta$ and \eqref{diff11} by $\dt\phi$ and add the resulting equations. 
After a cancellation of terms, we obtain the identity
\begin{align}\label{diff12} 
&\intQt|\dt\theta|^2\,+\,\frac 12\,\|A^\rho\gianni\theta(t)\|^2\,+\,\frac 12\,\|\dt\phi(t)\|^2+\pier{\int_0^t\juergen{\|}\dt\phi(s)\|^2_{V_B^\sigma} ds }\nonumber\\
&=\intQt u \, \dt\theta\,-\intQt(\ell(\phi_1)-\ell(\phi_2))\dt\phi_1\dt\theta\,-\intQt(F''(\phi_1)-F''(\phi_2))\dt\phi_1
\dt\phi\nonumber\\
&\hspace*{5mm}\pier{{}+\intQt (1- F''(\phi_2))}|\dt\phi|^2\,+\intQt (\ell'(\phi_1)-\ell'(\phi_2))\gianni{\theta_1\dt\phi_1}\dt\phi
\,+\intQt\ell'(\phi_2)\theta_1|\dt\phi|^2\nonumber\\
&\hspace*{5mm}+\intQt\dt\theta_1(\ell(\phi_1)-\ell(\phi_2))\dt\phi\,+\intQt\theta \, \ell'(\phi_2)\dt\phi_2\dt\phi
\,=:\,\sum_{j=1}^8 J_j\,,
\end{align}
with obvious notation. 
We estimate the terms on the right-hand side individually, using the Young and H\"older
inequalities, the embeddings from {\bf (A4)}, and the estimates \eqref{ssbound1}, \eqref{ssbound2}, and
\eqref{diff10} without further reference. Again, we denote by $C_i>0$, $i\in\enne$, constants that depend
only on $R$ and the data. 

Now let $\delta>0$ be arbitrary (to be chosen later). 
We obviously have 
\begin{align}\label{diff13}
|J_1|+|J_4|\,\le\,\delta\intQt|\dt\theta|^2\,+\,C_1\bigl(1+\delta^{-1}\bigr)\intQt|u|^2.
\end{align}
Moreover, \pier{it is clear that}
\begin{align}\label{diff14}
|J_2|\,&\le\,C_2\iot\|\dt\theta(s)\|\,\|\dt\phi_1(s)\|_{\Lx 4}\,\|\phi(s)\|_{\Lx 4}\,ds
\nonumber\\
&\le\,\delta\intQt|\dt\theta|^2\,+\,\frac{C_3}\delta \,\|\phi\|_{L^\infty(0,t;\VB\sigma)}^2
\iot\|\dt\phi_1(s)\|^2_{\VB \sigma}\,ds\nonumber\\
&\le\,\delta\intQt|\dt\theta|^2\,+\,\frac{C_4}\delta\intQt|u|^2\,.
\end{align}
Also, using \pier{the Lipschitz continuity of $F''$ in $[a_R, b_R]$ and} 
\eqref{diff10} once more, \pier{we infer} that
\begin{align}\label{diff15}
\gianni{|J_3|}\,&\le\,C_5\iot\|\phi(s)\|_{\Lx 4}\,\|\dt\phi_1(s)\|_{\Lx 4}\,\|\dt\phi(s)\|\,ds
\nonumber \\
&\le\,\pier{C_6}\,\|\phi\|_{L^\infty(0,t;\VB\sigma)}\,\|\phi_1\|_{H^1(0,t;\VB\sigma)}\,\|\phi\|_{H^1(0,t;H)}\nonumber\\
&\le\,\pier{C_7}\intQt|u|^2\,.
\end{align}
Similarly, \pier{in view of {\bf (F4)} we observe that}
\begin{align}\label{diff16}
|J_5|\,&\le\,\pier{C_8}\iot\|\phi(s)\|_{\Lx 4}\,\pier{\|\theta_1(s)\|_{\Lx 4} \,  \|\dt\phi_1(s)\|_{\Lx 4}}\,
\|\dt\phi(s)\|_{\Lx 4}\,ds\nonumber\\
&\le\,\delta\iot\|\dt\phi(s)\|^2_{\VB\sigma}\,ds\,+\,\frac{\pier{C_9}}\delta
\pier{\,\|\theta_1\|^2_{L^\infty(0,t;\VA\rho)}
\,\|\phi_1\|_{H^1(0,t;\VB\sigma)}^2
\,\|\phi\|_{L^\infty(0,t;\VB\sigma)}^2}
\nonumber\\
&\le\,\delta\pier{\iot\|\dt\phi(s)\|^2_{\VB\sigma}\,ds}\,+{}\frac{\pier{C_{10}}}\delta\intQt|u|^2\,.
\end{align}
We also have
\begin{align}\label{diff17}
|J_6|\,&\le\,\pier{C_{11}}\iot\|\theta_1(s)\|_{\Lx 4}\,\|\dt\phi(s)\|_{\Lx 4}\,\|\dt\phi(s)\|\,ds\nonumber\\ 
&\le\,\delta\iot\|\dt\phi(s)\|^2_{\VB\sigma}\,ds\,+\frac{\pier{C_{12}}}\delta\,\|\theta_1\|_{L^\infty(0,t;\VA\rho)}^2\iot
\|\dt\phi(s)\|^2\,ds\nonumber\\
&\le\,\delta\iot\|\dt\phi(s)\|^2_{\VB\sigma}\,ds\,+\,\frac{\pier{C_{13}}}\delta\intQt|u|^2\,.
\end{align}
Moreover, \pier{it turns out that}
\begin{align}\label{diff18}
|J_7|\,&\le\,\pier{C_{14}}\iot\|\dt\gianni{\theta_1}(s)\|\,\|\phi(s)\|_{\Lx 4}\,\|\dt\phi(s)\|_{\Lx 4}\, ds\nonumber\\
&\le\,\delta\iot\|\dt\phi(s)\|^2_{\VB\sigma}\,ds\,+\frac{\pier{C_{15}}}\delta\,\pier{\,\norma{\dt \theta_1}_{\Lt2H}^2}\,\|\phi\|_{L^\infty(0,t;\VB\sigma)}^2
\nonumber\\
&\le\,\delta\iot\|\dt\phi(s)\|^2_{\VB\sigma}\,ds\,+\,\frac{\pier{C_{16}}}\delta\intQt|u|^2\,.
\end{align}
Finally, \pier{we deduce that}
\begin{align}\label{diff19}
|J_8|\,&\le\,\pier{C_{17}}\iot\|\theta(s)\|\,\|\dt\phi_2(s)\|_{\Lx 4}\,\|\dt\phi(s)\|_{\Lx 4}\,ds\nonumber\\
&\le\,\delta\iot\|\dt\phi(s)\|^2_{\VB\sigma}\,ds\,+\,\frac{\pier{C_{18}}}\delta\,\|\theta\|_{L^\infty(0,t;H)}^2\iot   
\|\dt\phi_2(s)\|^2_{\VB\sigma}\,ds\nonumber\\
&\le\,\delta\iot\|\dt\phi(s)\|^2_{\VB\sigma}\,ds\,+\,\frac{\pier{C_{19}}}\delta\intQt|u|^2\,.
\end{align}
Summarizing the estimates \eqref{diff12}--\eqref{diff19}, and choosing $\delta>0$ small enough,
we have thus shown the estimate
\begin{align}\label{diff20}
\|\theta\|_{H^1(0,t;H)\cap L^\infty(0,t;\VA\rho)}^2\,+\,\|\phi\|_{W^{1,\infty}(0,t;H)\cap H^1(0,t;\VB\sigma)}^2
\,\le\,\pier{C_{20}}\intQt|u|^2\,.
\end{align} 
This concludes the proof of the assertion. 
\Edim


\section{Fr\'echet differentiability of ${\cal S}$}
\label{FRECHET}
\setcounter{equation}{0}

In this section, we aim to show the Fr\'echet differentiability of the control-to-state mapping ${\cal S}$
between suitable Banach spaces. 
To this end, we fix some $\,\bu\in{\cal U}_R\,$ and set $\,(\btheta,\bphi)={\cal S}(\bu)$. 
We then consider the linearized problem  
\begin{align}
\label{ls1}
&\dt\eta+\ell'(\bphi)\pier{\,\dt\bphi\,\xi}+\ell(\bphi)\dt\xi+A^{2\rho}\eta=h\quad\mbox{in }\,Q,\\[0.5mm]
\label{ls2}
&\dt\xi+B^{2\sigma}\xi+F''(\bphi)\xi=\ell'(\bphi)\,\btheta\,\xi+\ell(\bphi)\eta\quad\mbox{in }\,Q,\\[0.5mm]
\label{ls3}
&\eta(0)=\xi(0)=0\quad\mbox{in }\,\Omega.
\end{align}
The expectation is that if a Fr\'echet derivative $D{\cal S}(\bu)$ of ${\cal S}$ at $\bu$ exists, 
then, \gianni{for a given direction~$h$, it should
satisfy $D{\cal S}(\bu)[h]=(\eta,\xi)$, where} $(\eta,\xi)$ is the solution to \eqref{ls1}--\eqref{ls3}. 
We first show the following result.

\Bthm
\label{Wellposlin}
Suppose that the general \gianni{assumptions} {\bf (F1)}--{\bf (F4)}, {\bf (A1)}--{\bf (A4)} and {\bf (A6)} 
\gianni{as well as {\bf(GB)}} are fulfilled,
and let $\,\bu\in{\cal U}_R$ be arbitrary and $\,(\btheta,\bphi)={\cal S}(\bu)$. 
Then the linearized system
\eqref{ls1}--\eqref{ls3} has for every $\,h\in L^2(Q)\,$ a unique solution $\,(\eta,\pier{\xi})\,$ such that
\begin{align}
\label{regeta}
&\eta \in H^1(0,T;H)\cap L^\infty(0,T;\VA\rho)\cap L^2(0,T;\VA{2\rho}),\\
\label{regxi}
&\xi \in H^1(0,T;H)\cap L^\infty(0,T;\VB\sigma)\cap L^2(0,T;\VB {2\sigma}).
\end{align}
Moreover, the linear mapping $h\mapsto (\eta,\xi)$ is continuous as a mapping between the spaces
$L^2(Q)$ and $\bigl(H^1(0,T;H)\cap L^\infty(0,T;\VA\rho)\cap L^2(0,T;V_A^{2\rho})\bigr)\times \bigl((H^1(0,T;H)\cap L^\infty(0,T;\VB\sigma)\cap L^2(0,T;V_B^{2\sigma}\juerg{)}\bigr)$.
\Ethm

\Bdim                                    
We use a Faedo--Galerkin method. 
To this end, let (see \eqref{eigen}) $\{e_j\}_{j\in\enne}\,$ and 
$\,\{e_j'\}_{j\in\enne}\,$ be the orthonormalized eigenfunctions of $A$ and $B$, respectively. 
We define the $n$-dimensi\-onal spaces $V_n:={\rm span}\{e_1,\ldots,e_n\}$ and $V_n':={\rm span}\{e_1',\ldots,e_n'\}$ 
and search for every $n\in\enne$ functions of the form
$$\eta_n(x,t)=\sum_{j=1}^n a_j(t)e_j(x), \quad \xi_n(x,t)=\sum_{j=1}^n b_j(t)e_j'(x),$$
such that 
\begin{align}\label{gs1}
&(\dt\eta_n,v) + (A^\rho\eta_n,A^\rho v)+(\ell(\bphi)\dt\xi_n,v)\,=\,
-(\ell'(\bphi)\,\dt\bphi\,\xi_n,v)+(h,v)\nonumber\\[0.5mm]
&\qquad\mbox{for every $\,v\in V_n\,$ and a.e.\ in \,$(0,T)$,}\\[1mm]
\label{gs2}
&(\dt\xi_n,v)+(B^\sigma\xi_n,B^\sigma v) +(F''(\bphi)\xi_n,v)
\,=\,(\ell'(\bphi)\,\btheta\,\xi_n,v)+(\ell(\bphi)\eta_n,v)\nonumber\\[0.5mm]
&\qquad \mbox{for every $\,v\in V_n'\,$ and a.e.\ in }\,(0,T),\\[1mm]
\label{gs3}
&\eta_n(0)=\xi_n(0)=0.
\end{align}
We choose $v=e_k'$, $1\le k\le n$, in \eqref{gs2}, which, thanks to the orthogonality of the eigenfunctions, 
leads to $n$ explicit first-order ordinary differential equations with leading terms $\,\dt b_k$, $1\le k\le n$. 
Next,  we insert $v=e_k$, $1\le k\le n$, in \eqref{gs1}, 
and we substitute the explicit expressions for $\,\dt b_k$, $1\le k\le n$, in the terms $(\ell(\bphi)\dt\xi_n,e_k)$ for $k=1,...,n$. 
By doing this, we obtain from \eqref{gs1}--\eqref{gs3} a standard initial value problem 
for a linear system of $2n$ ordinary differential equations in the unknowns $a_1,\ldots,a_n,b_1,\ldots,b_n$. 
Since all of
the occurring coefficient functions belong to $L^2(0,T)$, it follows from Carath\'eodory's theorem the
existence of a unique solution $(a_1,\ldots,a_n,b_1,\ldots,b_n)\in H^1(0,T;\erre^{2n})$ which specifies the
unique solution $(\eta_n,\xi_n)\in (H^1(0,T;V_n)\times H^1(0,T;V_n')\juerg{)}$ to \eqref{gs1}--\eqref{gs3}.

In the following, we derive some a priori estimates for the Galerkin approximations. 
In this process, we denote by $C_i>0$, $i\in\enne$, constants that may depend on $R$ and the data of the state
system, but not on $n\in\enne$. 
\gianni{We recall that Remark~\ref{FromGB} applies to~$\bphi$.}
To begin with, we insert $v=\etan$ in \eqref{gs1} and $v=\dt\xin$ in \eqref{gs2},
and add the results, which leads to a cancellation of terms. 
Then, we integrate over time and add
to both sides of the resulting identity the same term \gianni{$\frac 12\,\norma{\xin(t)}^2=\intQt \xin\dt\xin$}. 
We then obtain the equation
\begin{align}
\label{gs4}
&\frac 12 \,\|\etan(t)\|^2\,+\,\frac 12\,\|\xin(t)\|^2_{\VB\sigma}\,+\intQt|A^\rho\etan|^2\,+\intQt|\dt\xin|^2
\nonumber\\
&=\,\intQt h\etan\,-\intQt \ell'(\bphi)\,\dt\bphi\,\xin\etan\,-\intQt F''(\bphi)\xin\,\dt\xin\nonumber\\
&\hspace*{5mm}+\intQt\ell'(\bphi)\,\btheta\,\xin\,\dt\xin\,+\intQt\xin\,\dt\xin\,=:\,\sum_{j=1}^{\gianni5} L_j,
\end{align}
with obvious meaning. 
Let $\delta>0$ be arbitrary (to be specified later). 
At first, it is readily seen that
\begin{align}\label{gs5}
|L_1|+|L_3|+|L_5|\,\le\,\frac 12\intQt(|h|^2+|\etan|^2)\,+\,\delta\intQt|\dt\xin|^2\,+\,\frac{C_1}\delta\intQt|\xin|^2\,.
\end{align}
Moreover, \pier{we observe that}
\begin{align}\label{gs6}
|L_2|\,&\le\,C_2\iot\|\dt\bphi(s)\|_{\Lx 4}\,\|\xin(s)\|_{\Lx 4}\,\|\etan(s)\|\,ds\nonumber\\
&\le\,C_3\intQt|\etan|^2\,+\,C_4\iot\|\dt\bphi(s)\|^2_{\VB\sigma}\,\|\xin(s)\|^2_{\VB\sigma}\,ds\,.
\end{align}
Finally, we have the estimate
\begin{align}\label{gs7}
|L_4|\,&\le\,C_5\iot\|\btheta(s)\|_{\Lx 4}\,\|\xin(s)\|_{\Lx 4}\,\|\dt\xin(s)\|\,ds\nonumber\\
&\le\,\delta\intQt|\dt\xin|^2\,+\,\frac{C_6}\delta\iot\|\btheta(s)\|^2_{\VA\rho}\,\|\xin \pier{(s)}\|^2_{\VB\sigma}\,ds\,.
\end{align}
Now observe that the mapping $\,\,s\mapsto \|\dt\bphi(s)\|^2_{\VB\sigma}+\|\btheta(s)\|^2_{\VA\rho}\,\,$
is known to belong to $L^1(0,T)$. 
Hence, combining \eqref{gs4}--\eqref{gs6}, and
choosing $\delta>0$ small enough, we obtain from Gronwall's lemma
the estimate
\begin{equation}\label{gs8}
\|\etan\|_{L^\infty(0,T;H)\cap L^2(0,T;\VA\rho)}\,+\,\|\xin\|_{H^1(0,T;H)\cap L^\infty(0,T;\VB\sigma)}\,\le\,
C_7\|h\|_{L^2(0,T;H)}\,.
\end{equation}

From \eqref{gs8} we can draw some consequences. 
Namely, invoking the general bounds \eqref{ssbound1} and \eqref{ssbound2},
\gianni{\juerg{as well as} the embeddings given by~{\bf(A4)}}, we can easily verify that 
\begin{align*}
&\|\ell'(\bphi)\,\dt\bphi\,\xin+\ell(\bphi)\,\dt\xin\|_{L^2(0,T;H)}\,\le\,C_8\|h\|_{L^2(0,T;H)},\\[0.5mm]
&\|\ell'(\bphi)\,\btheta\,\xin+\ell(\bphi)\,\etan-F''(\bphi)\,\xin\|_{L^2(0,T;H)}\,\le\,C_9
\|h\|_{L^2(0,T;H)}.
\end{align*}
But then we may insert first $v=\dt\etan$ and then $v=A^{2\rho}\etan$ in \eqref{gs1} to conclude that 
\begin{equation}\label{gs9}
\|\etan\|_{H^1(0,T;H)\cap L^\infty(0,T;\VA\rho)\cap L^2(0,T;\VA {2\rho})}
\,\le\,C_{10}\|h\|_{L^2(0,T;H)}.
\end{equation}
Likewise, by inserting $v=B^{2\sigma}\xin$ in \eqref{gs2}, we find that 
\begin{equation}\label{gs10}
\|\xin\|_{H^1(0,T;H)\cap L^\infty(0,T;\VB\sigma)\cap L^2(0,T;\VB{2\sigma})}\,\le\,
C_{11}\|h\|_{L^2(0,T;H)}.
\end{equation}
Hence, there is a pair $(\eta,\xi)$ such that (first only for a subsequence, but, by the uniqueness of the limit,
eventually for
the entire sequence) we have the convergence properties
\begin{align*} 
&\etan\to\eta\quad\mbox{weakly-star in }\,H^1(0,T;H)\cap L^\infty(0,T;\VA\rho)\cap L^2(0,T;\VA {2\rho}),\nonumber\\[0.5mm]
&\xin\to\xi\quad\mbox{weakly-star in }\,H^1(0,T;H)\cap L^\infty(0,T;\VB\sigma)\cap L^2(0,T;\VA {2\sigma}). 
\end{align*}                       
It is then a standard matter (which needs no repetition here) to show that $(\eta,\xi)$ is a strong solution
to the linearized system \eqref{ls1}--\eqref{ls3}, and the validity of the assertion concerning the continuity of  
the mapping $\,h\mapsto (\eta,\xi)\,$ follows from \eqref{gs9} and \eqref{gs10} by using the semicontinuity properties
of norms.

It remains to show the uniqueness of the solution. To this end, let $(\eta_i,\xi_i)$, $i=1,2$, be two solutions 
enjoying the regularity properties \eqref{regeta} and \eqref{regxi}, and let $\eta:=\eta_1-\eta_2$, $\xi:=\xi_1-\xi_2$. 
Then $(\eta,\xi)$ is a strong solution to the system \eqref{ls1}--\eqref{ls3} with $h=0$. 
Repeating the a priori estimates leading to \eqref{gs8} for the continuous problem, we obtain an estimate for 
$(\eta,\xi)$ which resembles \eqref{gs8}, but this time with $h=0$ on the right-hand side. 
Thus, $\eta=\xi=0$. 
This concludes the proof of the assertion.
\Edim              

We are now ready to prove the Fr\'echet differentiability of the control-to-state operator. 
To be able to perform this analysis, we need a slightly stronger embedding condition than that of assumption {\bf (A4)}. 
We have to postulate:

\vspace{2mm}\noindent                          
{\bf (A8)} \,\,\,The embeddings $\,V_A^{2\rho}\subset {\Lx 6}\,$ and $\,V_B^\sigma\subset \Lx 6$ are continuous.

\vspace{2mm}
\Brem 
\label{RemA8}
The second condition is more restrictive than the first one. 
Indeed, if $A=B=-\Delta$ with zero Dirichlet or Neumann conditions, 
then $V_B^\sigma  \subset H^{2\sigma}(\Omega)\subset \Lx 6$ if $\sigma\ge 1/2$, by \eqref{sobslo3}. 
On the other hand, $V_A^{2\rho}\subset H^{4\rho}(\Omega)\subset \Lx 6$ provided that \gianni{$\rho\ge 1/4$}. 
Notice that \gianni{$V_A^\rho\subset \Lx 4$ if $\rho\ge 3/8$ (see also Remark~\ref{RemA4})}, 
so that in this case the postulate 
for $A$ in {\bf (A8)} is not more restrictive than that required in~{\bf (A4)}.
\Erem
  
With these preparations, the road is paved for the proof of Fr\'echet differentiability.

\vspace{1mm}
\Bthm
\label{Frechet}
Suppose that the conditions {\bf (F1)}--{\bf (F4)}, {\bf (A1)}--{\bf (A4)}, {\bf (A6)}--{\bf (A8)},
and {\bf (GB)} are fulfilled. 
Then the control-to-state operator ${\cal S}$ is Fr\'echet differentiable on 
${\cal U}_R$ as a mapping from $L^\infty(Q)$ into the Banach space
\begin{equation}\label{defY}
{\cal Y}:=\bigl(H^1(0,T;V_A^{-\rho})\cap C^0([0,T];H)\cap L^2(0,T;\VA\rho)\bigr)\times\bigl(H^1(0,T;H)\cap L^\infty(0,T;\VB\sigma)\bigr).
\end{equation}
Moreover, if $\,\bu\in{\cal U}_R\,$ and $(\btheta,\bphi)={\cal S}(\bu)$, then the Fr\'echet derivative 
\juergen{$\,D{\cal S}(\bu)\in{\cal L}(L^\infty(Q),{\cal Y})$ of $\,{\cal S}\,$ at $\,\bu$, applied to $\,h\in L^\infty(Q)$,} 
satisfies $\,D{\cal S}(\bu)[h]=(\eta,\xi)$, where $\,(\eta,\xi)\,$ is the unique solution of \eqref{ls1}--\eqref{ls3}.
\Ethm 
                  
\Bdim
Let $\bu\in{\cal U}_R$ be arbitrary and $(\btheta,\bphi)={\cal S}(\bu)$. 
Since ${\cal U}_R$ is open, there is some
$\Lambda>0$ such that $\bu+h\in{\cal U}_R$ whenever $\|h\|_{L^\infty(Q)}\le\Lambda$. 
In the following, we only consider such perturbations~$\,h$. 
For any such $\,h$, we set $\,(\thetah,\phih)={\cal S}(\bu+h)$, and
we denote by $(\etah,\xih)$ the unique solution to the linearized system \eqref{ls1}--\eqref{ls3} associated with~$\,h$. 
Moreover, we~put
$$
y^h:=\thetah-\btheta-\etah,\quad z^h:=\phih-\bphi-\xih.
$$
Since the linear mapping $h\mapsto (\etah,\xih)$ is by Theorem~\ref{Wellposlin} continuous from $\Liq$ into ${\cal Y}$, 
it suffices to show that there is a mapping $\,Z:(0,+\infty)\to (0,+\infty)$ such that
\begin{equation}\label{f1}
\|(y^h,z^h)\|_{\cal Y}\,\le\,Z(\|h\|_{\Liq})\,\,\,\,\mbox{ and }\,\,\,\,\lim_{s\searrow0}\frac{Z(s)}s\,=0.
\end{equation}

To begin with, note that $(y^h,z^h)$ satisfies the regularity properties
\gianni{(see also \eqref{moreregphi})}
\begin{align}\label{regyh}
&y^h \in H^1(0,T;H)\cap L^\infty(0,T;\VA\rho)\cap L^2(0,T;\VA{2\rho}),\\[0.5mm]
\label{regzh}
&z^h\in H^1(0,T;H)\cap L^\infty(0,T;\VB\sigma)\cap L^2(0,T;\VB{2\sigma}).
\end{align}
Moreover, $(\btheta,\bphi)$ and $(\thetah,\phih)$ satisfy the global  estimates
\eqref{ssbound1} and \eqref{ssbound2}, and from \eqref{ssbound3} we have for all $t\in (0,T]$ the estimate
\begin{align}\label{f2}
\|\thetah-\btheta\|_{H^1(0,t;H)\cap L^\infty(0,t;\VA\rho)}\,+\,\|\phih-\bphi\|_{W^{1,\infty}(0,t;H)
\cap H^1(0,t;\VB\sigma)}\,\le\,K_2\,\|h\|_{L^2(0,t;\gianni H)}\,.
\end{align}                         
In the following, we denote by $C>0$ constants that may depend on $R$ and the data of the state system, but not on the 
special choice of $h\in\Liq$ with $\bu+h\in{\cal U}_R$. 
Observe that the meaning of $\,C\,$ may change from line to line within formulas. 

At this point, we observe that Taylor's theorem with integral remainder shows that we have almost everywhere in $\,Q\,$
the identities 
\begin{align}\label{Taylor1}
\ell(\phih)\,&=\,\ell(\bphi)+\ell'(\bphi)(\phih-\bphi)+(\phih-\bphi)^2\,R_1^h,\\[0.5mm]
\label{Taylor2}
F'(\phih)\,&=\,F'(\bphi)+F''(\bphi)(\phih-\bphi)+(\phih-\bphi)^2\,R_2^h,
\end{align}   
with the remainders
\begin{align*}
&R_1^h=\int_0^1 (1-s)\ell''(\bphi+s(\phih-\bphi))\,ds, \quad 
R_2^h=\int_0^1 (1-s)F''' (\bphi+s(\phih-\bphi))\,ds\,.
\end{align*}
By \pier{{\bf (F2)}--{\bf (F4)}, {\bf (GB)} and the boundedness of $F''' $ in $[a_R, b_R]$,}  we have 
\begin{equation}\label{Taylor3}
\|R_1^h\|_{\Liq}\,+\,\|R_2^h\|_{\Liq}\,\le\,C.
\end{equation}
Now observe that $\,y^h\,$ and $\,z^h\,$ are strong solutions to the system
\begin{align}\label{f3}
&\dt y^h+A^{2\rho}y^h \,=\,Q_1^h\quad\mbox{in }\,Q,\\[0.5mm]
\label{f4}
&\dt z^h+B^{2\sigma}z^h\,=\,Q_2^h\quad\mbox{in }\,Q,\\[0.5mm]
\label{f5}
&y^h(0)=z^h(0)=0\quad\mbox{in }\,\Omega, 
\end{align}
where simple algebraic manipulations using \eqref{Taylor1} and \eqref{Taylor2} show that
\begin{align}\label{Q1h}
Q_1^h\,&=\,-(\ell(\phih)-\ell(\bphi))(\dt\phih-\dt\bphi)-\ell(\bphi)\dt z^h-\ell'(\bphi)z^h\,\dt\bphi
-R_1^h\,(\phih-\bphi)^2\,\dt\bphi,
\\[0.5mm]
\label{Q2h}
Q_2^h\,&=\,(\ell(\phih)-\ell(\bphi))(\thetah-\btheta)+\ell(\bphi)y^h+\btheta\,\ell'(\bphi)z^h
+\btheta\,R_1^h\,(\phih-\bphi)^2-F''(\bphi)z^h\nonumber\\
&\quad\,\,-R_2^h\,(\phih-\bphi)^2\,.
\end{align}
Now we test \eqref{f3} by $\,y^h\,$ and \eqref{f4} by $\,\dt z^h$, add the results (whereby two terms cancel),
and add the same \gianni{term $\,\,\frac12\,\norma{z^h(t)}^2=\intQt z^h\dt z^h\,\,$} to both sides of the resulting identity. 
\gianni{Since the terms involving the product $y^h\dt z^h$ cancel out, we} obtain that
\begin{align}
\label{f6}
&\frac 12\,\|y^h(t)\|^2\,+\,\frac 12\,\|z^h(t)\|^2_{\VB\sigma}\,+\intQt|\dt z^h|^2\,+\intQt|A^\rho y^h|^2\nonumber\\[1mm]
&=\,-\intQt\!(\ell(\phih)-\ell(\bphi))(\dt\phih-\dt\bphi)\,y^h-\intQt\! \ell'(\bphi)\, \pier{\dt\bphi\,z^h}\,y^h \nonumber\\[1mm]
&\quad\,\, 
-\intQt\! R_1^h\,(\phih-\bphi)^2\,\dt\bphi\,y^h\nonumber
+\intQt(\ell(\phih)-\ell(\bphi))(\thetah-\btheta)\,\dt z^h\\[1mm]
&\quad\,\, +\intQt\btheta\,\ell'(\bphi)z^h\,\dt z^h
\,+\intQt\btheta\,R_1^h\,(\phih-\bphi)^2\,\dt z^h\nonumber\\[1mm]
&\quad\,\,+\intQt (1-F''(\bphi))z^h\,\dt z^h\,-\intQt R_2^h\,(\phih-\bphi)^2\,\dt z^h\,\,=:\,\,\sum_{j=1}^8M_j,
\end{align}   
with obvious meaning. 
Let $\delta>0$ be arbitrary (to be specified later). 
We estimate the terms on the right-hand side individually, using the Young and H\"older inequalities, 
the global bounds \eqref{ssbound1}, \eqref{ssbound2} and \eqref{Taylor3}, 
the stability estimate \eqref{f2}, as well as the embedding conditions {\bf (A4)} and~{\bf (A8)},
repeatedly without further reference. 
Here, for the sake of brevity,
we often omit the argument $\,s\,$ of the involved functions. 
At first, we have
\begin{align}
\label{f7}
|M_1|\,&\le\,C\iot\|\phih-\bphi\|_{\Lx 4}\,\|\dt\phih-\dt\bphi\|_{\Lx 4}\,\|y^h\|\,ds\nonumber\\
&\le\,C\,\|\phih-\bphi\|_{L^\infty(0,t;\VB\sigma)}^2\,\|\phih-\bphi\|_{H^1(0,t;\VB\sigma)}^2\,+
\iot\|y^h\|^2\,ds
\nonumber\\
&\le\,C\,\|h\|^4_{L^2(0,t;\gianni H)}\,+\iot\|y^h\|^2\,ds\,.
\end{align}
Moreover, \pier{we see that}
\begin{align}\label{f8}
|M_2|\,&\le\,C \iot \pier{ \|\dt\bphi\|_{\Lx 4}\,\|z^h\|_{\Lx 4}}\,\|y^h\|\,ds\nonumber\\
&\le\, \pier{C\iot\|\dt\bphi\|_{\VB\sigma} \left( \|y^h\|^2 + \|z^h\|^2_{\VB\sigma}\right)}ds, 
\end{align}
as well as 
\begin{align}\label{f9}
|M_3|\,&\le\,C\intQt \|\phih-\bphi\|^2_{\Lx 6}\,\|\dt\bphi\|_{\Lx 6}\,\|y^h\|\,ds\nonumber\\
&\le\,C\,\|\phih-\bphi\|^4_{L^\infty(0,T;\VB\sigma)}\,\|\bphi\|^2_{H^1(0,T;\VB\sigma)}\,+\,\iot\|y^h\|^2\,ds
\nonumber\\
&\le\,C\,\|h\|^4_{L^2(0,t;H)}+\,\iot\|y^h\|^2\,ds\,.
\end{align}
Also, \pier{it follows that}
\begin{align}\label{f10}
|M_4|\,&\le\,C\iot \|\phih-\bphi\|_{\Lx 4}\,\|\thetah-\btheta\|_{\Lx 4}\,\|\dt z^h\|\,ds \nonumber\\
&\le\,\delta\iot\|\dt z^h\|^2\,ds\,+\,\frac C\delta\,\|\phih-\bphi\|_{L^\infty(0,t;\VB\sigma)}^2\,
\|\thetah-\btheta\|^2_{L^\infty(0,t;\VA\rho)}\nonumber\\
&\le\,\delta\iot\|\dt z^h\|^2\,ds\,+\,\frac C\delta\,\|h\|_{L^2(0,t;H)}^4\,,
\end{align}
and \pier{that}
\begin{align}\label{f11}
|M_5|\,&\le\,C\iot \|\btheta\|_{\Lx 4}\,\|z^h\|_{\Lx 4}\,\|\dt z^h\|\,ds\,\le\,\delta
\iot\|\dt z^h\|^2\,ds\,+\,\frac C\delta\iot\|z^h\|_{\VB\sigma}^2\,ds\,,
\end{align}
as well as
\begin{align}\label{f12}
|M_6|\,&\le\,C\iot \|\btheta\|_{\Lx 6}\,\|\phih-\bphi\|^2_{\Lx 6}\,\|\dt z^h\|\,ds\nonumber\\
&\le\,\delta\iot\|\dt z^h\|^2\,ds\,+\,\frac C\delta\,\|\phih-\bphi\|^4_{L^\infty(0,t;\VB\sigma)}\nonumber\\
&\le\,\delta\iot\|\dt z^h\|^2\,ds\,+\,\frac C\delta\,\|h\|^4_{L^2(0,t;H)}\,.
\end{align}
Finally, \pier{we infer that}
\begin{align}\label{f13}
|M_7|\,\le\,\delta\iot\|\dt z^h\|^2\,ds \,+\,\frac C\delta\iot\|z^h\|^2\,ds
\end{align}
and 
\begin{align}\label{f14}
|M_8|\,&\le\, C\iot\|\phih-\bphi\|_{\Lx 4}^2\,\|\dt z^h\|\,ds
\,\le\,\delta\,\iot\|\dt z^h\|^2\,ds\,+\,\frac C\delta\,\|h\|^4_{L^2(0,t;H\pier )}\,.
\end{align}

At this point, \gianni{we observe that the map
\pier{$\,\,s\mapsto\|\dt\bphi(s)\|_{\VB\sigma}\,$ belongs to $L^2(0,T)$}. 
Thus, choosing $\delta>0$ small enough and combining the estimates \accorpa{f6}{f14}}, we conclude that
\begin{align}
\label{f15}
\|y^h\|_{L^\infty(0,T;H)\cap L^2(0,T;\VA\rho)}\,+\,\|z^h\|_{H^1(0,T;\gianni H)\cap L^\infty(0,T;\VB\sigma)}\,\le\,
C\,\|h\|^2_{L^2(0,T;H)}\,.
\end{align}  
With this estimate shown, it is a simple comparison argument in \eqref{f3} (which we may leave to the reader)
to verify that also
\begin{equation*}
\|y^h\|_{H^1(0,T;\VA{-\rho})}\,\le\,C\,\|h\|^2_{L^2(0,T;H)}\,.
\end{equation*}  
Now observe that $H^1(0,T;\VA{-\rho})\cap L^2(0,T;\VA\rho)$ is continuously embedded in $C^0([0,T];H)$, \pier{so that 
\eqref{f1} is satisfied} with a function of the form $Z(s)=\hat Cs^2$, \pier{for} a sufficiently large $\hat C>0$. 
This concludes the proof of the assertion.
\Edim

As an immediate consequence of Theorem~\ref{Frechet}, we now deduce a first necessary optimality condition for the optimal
control problem~{\bf (CP)}.

\Bcor
\label{Badoc}
Suppose that the conditions {\bf (F1)}--{\bf (F4)}, {\bf (A1)}--{\bf (A8)}, and {\bf (GB)} hold true. 
Moreover, let $\,\bu\in\uad\,$ be an optimal control for the problem {\bf (CP)} and $\,(\btheta,\bphi)={\cal S}(\bu)$. 
Then, there holds the variational inequality
\begin{align}\label{vug1}
&\beta_1\iO(\bphi(T)-\phi_\Omega)\,\xi(T)\,+\,\beta_2\iint_Q(\bphi-\phi_Q)\,\xi\,+\,\beta_3\iO(\btheta(T)-\theta_\Omega)
\,\eta(T)\nonumber\\
&+\,\beta_4\iint_Q(\btheta-\theta_Q)\,\eta\,+\,\beta_5\iint_Q \bu\,(u-\bu)\,\ge\,0\qquad\forall\,u\in\uad,
\end{align}
where $(\eta,\xi)$ is the unique solution to the linearized system \eqref{ls1}--\eqref{ls3} associated with $h=u-\bu$. 
\Ecor

\Bdim
By virtue of the quadratic form \juerg{of} $\,{\cal J}\,$ and Theorem~\ref{Frechet}, the reduced cost functional
$\,\widetilde {\cal J}(u):={\cal J}({\cal S}(u),u)\,$ is Fr\'echet differentiable on ${\cal U}_R$. Since $\uad$ is
convex, we must have $D\widetilde{\cal J}(\bu)[u-\bu]\ge 0$ for all $u\in\uad$. The result then  follows in a
standard manner from
the chain rule and Theorem~\ref{Frechet}.
\Edim


\section{The optimal control problem}
\label{OPTIMACP}
\setcounter{equation}{0}

In this section, we investigate the optimal control problem {\bf (CP)}.


\subsection{Existence of optimal controls}
\label{EXISTOC}

We begin our analysis of {\bf (CP)} with an existence result.

\Bthm
\label{Existoc}
Suppose that {\bf (F1)}--{\bf (F4)} and {\bf (A1)}--{\bf (A4)} are fulfilled. 
Then {\bf (CP)} has a solution.
\Ethm

\Bdim
We pick a minimizing sequence $\{u_n\}\subset \uad$ and set $(\thetan,\phin):={\cal S}(\un)$, for all $n\in\enne$.
\gianni{We fix $R>0$ such that $\,\uad\subset {\cal U}_R$ and account for~\eqref{ssbound1}.
Hence,} invoking standard compactness results (cf., e.g., \cite[Sect.~8, Cor.~4]{Simon} for 
\gianni{the strong compactness}), 
we may assume that there are $u\in\uad$ and $(\theta,\phi)$ such that\pier{, at least
for a subsequence,}
\begin{align}
\label{convu}
\un\to u\quad&\mbox{weakly-star in }\,\Liq,\\[0.5mm]
\label{convt}
\thetan\to\theta\quad&\mbox{weakly-star in }\,H^1(0,T;H)\cap L^\infty(0,T;\VA\rho)\cap L^2(0,T;\VA{2\rho}),
\nonumber\\
&\mbox{strongly in $\,C^0([0,T];H)\,$ and pointwise a.e.\ in }\,Q,\\[0.5mm]
\label{convp}
\phin\to\phi\quad&\mbox{weakly-star in }\,\gianni{W^{1,\infty}(0,T;H)\cap H^1(0,T;\VB\sigma)},
\nonumber\\
&\mbox{strongly in $\,C^0([0,T];H)\,$ and pointwise a.e.\ in }\,Q.
\end{align}

We now show that $(\theta,\phi)={\cal S}(u)\,$ which implies that the pair $((\theta,\phi),u)$ is admissible for~{\bf (CP)}.
Once this is proved, the lower semicontinuity of norms shows that $((\theta,\phi),u)$ is an optimal pair.
  
At first, note that obviously $\theta(0)=\theta_0$ and $\phi(0)=\phi_0$. 
In additon, by Lipschitz continuity, 
$F_2'(\phin)\to F_2'(\phi)$\, and \,$\ell(\phin)\to\ell(\phi)$, both strongly in $C^0([0,T];H)$. 
Since $\{\ell(\phin)\thetan\}$ is easily seen to be bounded in $L^2(Q)$, the latter entails that   
$\ell(\phin)\thetan\to\ell(\phi)\theta\,$ weakly in $L^2(0,T;H)$. 

Now, we write the time-integrated version of \eqref{ssvar1}, written for $u=\un$ and $(\theta,\phi)=(\thetan,\phin)$,
with test functions $v\in L^2(0,T;\VA\rho)$. 
Taking the limit as $n\to\infty$, we find that $(\theta,\phi)$ 
satisfies the time-integrated version of \eqref{ssvar1}, which is equivalent to \eqref{ssvar1} \gianni{itself}. 
It remains to show the validity of \eqref{ssvar2}. 
To this end, we use the semicontinuity of $F_1$ and \eqref{convp}, 
which yield that $\,0\le F_1(\phi)\le\liminf_{n\to\infty} \,F_1(\phin)\,$ a.e.\ in $\,Q$.
Then, owing to Fatou's lemma and to~\eqref{ssbound1},
$$
0\,\le \juerg{\iint_Q} F_1(\phi)\,\le\,\liminf_{n\to\infty}\juerg{\iint_Q} F_1(\phin)\,\le\,K_1.
$$
In particular, $\,F_1(\phi)\in L^1(Q)$. 
Now, the quadratic form \gianni{$\,v\mapsto\ioT \norma{B^\sigma v(t)}^2\,dt\,$} is lower semi\-continuous on $L^2(0,T;\VB\sigma)$. 
\gianni{Thus, starting from \eqref{intssvar2} written for $(\theta_n,\phi_n)$, 
we can deduce, for every $\,v\in L^2(0,T;\VB\sigma)$, the following chain}:
\begin{align}
\label{exivug}
&\ioT(B^\sigma\phi(t),B^\sigma(\phi(t)-v(t)))\,dt\,+\juerg{\iint_Q} F_1(\phi)\nonumber\\
&\le\,\liminf_{n\to\infty} \,\Bigl(\ioT(B^\sigma\phin(t),B^\sigma(\phin(t)-v(t))\,dt\,+\juerg{\iint_Q} F_1(\phin)\Bigr)\nonumber\\
&\le\,\liminf_{n\to\infty} \,\Bigl(\juerg{\iint_Q}(\ell(\phin)\thetan-\dt\phin-F_2'(\phin))(\phin-v)
\,+\juerg{\iint_Q} F_1(v)\Bigr)\nonumber\\
&\gianni{{}={}}\,\juerg{\iint_Q}(\ell(\phi)\theta-\dt\phi-F_2'(\phi))(\phi-v)\,+\juerg{\iint_Q} F_1(v)\,.
\end{align}
In other words, $(\theta,\phi)$ satisfies \gianni{\eqref{intssvar2}}, which is equivalent to \eqref{ssvar2}. 
This concludes the proof of the assertion.
\Edim


\subsection{Necessary optimality conditions}
\label{NECOPT}

We now turn our interest to the derivation of first-order necessary optimality conditions. 
To this end, we assume that $\bu\in\uad$ is an optimal control with associated state $\,(\btheta,\bphi)$, 
and we assume that all of the general assumptions {\bf (F1)}--{\bf (F4)}, {\bf (A1)}--{\bf (A8)}, and {\bf (GB)} are satisfied. 
Hence, in particular, the double obstacle potential $F_{\rm 2obs}$ is excluded from the consideration. 
We aim at eliminating the expressions involving $(\eta,\xi)$ from the variational inequality \eqref{vug1} by means of the 
adjoint state variables. 
The adjoint system formally reads:
\begin{align}
\label{as1}
&-\dt q-\ell(\bphi)p+A^{2\rho}q\,=\,\beta_4(\btheta-\theta_Q)\quad\mbox{in }\,Q,\\[0.5mm]
\label{as2}
&-\dt p-\ell(\bphi)\dt q+B^{2\sigma}p +F''(\bphi)p-\ell'(\bphi)\btheta\,p\,=\,\beta_2(\bphi-\phi_Q)\quad
\mbox{in }\,Q,\\[0.5mm]
\label{as3}
&q(T)=\beta_3(\btheta(T)-\theta_\Omega), \quad p(T)=\beta_1(\bphi(T)-\phi_\Omega)-
\beta_3\ell(\bphi(T))(\btheta(T)- \theta_\Omega)\quad
\mbox{in }\,\Omega.
\end{align}
Owing to the low regularity of the final \gianni{data appearing in} \eqref{as3}, we cannot expect to obtain a strong solution to this system. 
Indeed, it turns out that 
\eqref{as2} is meaningful only in its weak form
\begin{align}
\label{as2w}
&\langle -\dt p(t),v\rangle_{\VB\sigma} \,-\,(\ell(\bphi(t))\dt q(t),v)+(B^\sigma p(t),B^\sigma v)
+(F''(\bphi(t))p(t),v)\nonumber\\
&\quad {}-(\ell'(\bphi(t))\btheta(t)p(t),v)\,=\,\beta_2(\bphi(t)-\phi_Q(t),v)
\non\\
&\quad\mbox{for all $\,v\in\VB\sigma\,$ and a.e.\ $t\in(0,T)$}.
\end{align}

Another point is that, in order to derive a priori bounds, one would like to test \eqref{as2} by $\,p\,$ and \eqref{as1}
by $\,-\dt q$, which makes it necessary to assume that the associated final datum belongs to $\VA \rho$. 
We thus postulate:

\vspace{2mm}\noindent
{\bf (A9)} \,\,\,It holds $\,\beta_3\,\theta_\Omega\in\VA\rho$.

\vspace{2mm}
\noindent 
This condition is satisfied if  $\beta_3=0$ or $\theta_\Omega\in \VA\rho$. 
In the first case, there is no endpoint tracking of the temperature in the cost functional, 
while in the second the regularity of the target function coincides with that of the associated state (which makes sense). 
We have the following well-posedness result:

\Bthm
\label{Wellposadj}
Let the assumptions {\bf (F1)}--{\bf (F4)}, {\bf (A1)}--{\bf (A9)}, and {\bf (GB)} be fulfilled, and let $\bu\in\uad$ be given
with associated state $(\btheta,\bphi)={\cal S}(\bu)$. Then the adjoint problem \eqref{as1}, \eqref{as2w}, 
\eqref{as3} has a unique solution $(p,q)$ such that
\begin{align}
\label{regq}
&q\in H^1(0,T;H)\cap C^0([0,T];\VA\rho)\cap L^2(0,T;\VA{2\rho}),\\[0.5mm]
\label{regp}
&p\in H^1(0,T;\VB {-\sigma})\cap C^0([0,T];H)\cap L^2(0,T;\VB\sigma).
\end{align}
\Ethm

\Bdim
As in the proof of Theorem~\ref{Wellposlin}, we use a Faedo--Galerkin technique with the eigenfunctions of the operators $A$ and~$B$.
\gianni{With the notations used there, we look} for functions of the form 
$$
q_n(x,t)=\sum_{j=1}^n a_j(t)e_j(x),\quad\pen(x,t)=\sum_{j=1}^n b_j(t)e_j'(x),
$$
satisfying the system
\begin{align}
\label{disc1}
&-(\dt\qn(t),v)-(\ell(\bphi(t))\pen(t),v)+(A^\rho\qn(t),A^\rho v)\,=\,(g_4(t),v)\nonumber\\[0.5mm]
&\qquad\mbox{for all $\,v\in \gianni{V_n}\,$ and a.e.\ $\,t\in(0,T)$,}\\[1mm]
\label{disc2}
&-(\dt\pen(t),v)-(\ell(\bphi(t))\dt\qn(t),v)+(B^\sigma\pen(t),B^\sigma v)+(F''(\bphi(t))\pen(t),v)\nonumber\\[0.5mm]
&\quad-(\ell'(\bphi(t))\btheta(t)\pen(t),v)\,=\,(g_2(t),v) \quad\mbox{for all\,$v\in \gianni{V_n'}\,$ and
a.e.\ $\,t\in(0,T)$,}\\[1mm]
\label{disc3}
&(\qn(T),v)=(g_3,v) \quad\forall\,v\in\gianni{V_n}\,,
\quad 
(\pen(T),v)=(g_1,v)-(\ell(\bphi(T))g_3,v)\quad\forall\,v\in\gianni{V_n'}\,,
\end{align} 
where we have set
\begin{align}
\label{gs}
&g_1=\beta_1(\bphi(T)-\phi_\Omega),\enskip g_2=\beta_2(\bphi-\phi_Q),\enskip g_3=\beta_3(\btheta(T)-\theta_\Omega),
\enskip g_4=\beta_4(\btheta-\theta_Q).
\end{align}
Using an analogous argument as in the proof of Theorem~\ref{Wellposlin}, we can infer that the system \eqref{disc1}--\eqref{disc3}
enjoys a unique solution pair $(\qn,\pen)\in (H^1(0,T;V_n)\times H^1(0,T;V_n'))$. 

We now derive a priori estimates for
the approximations $(\qn,\pen)$, where we denote by $C_i$, $i\in\enne$, constants that may depend on $R$ and the data, but not on $n\in\enne$. 
To begin with, we insert $\,v=-\dt\qn(t)\,$ in \eqref{disc1} and $\,v=\pen(t)\,$ in \eqref{disc2}, 
add the results, and integrate over $(t,T)$ where $t\in [0,T)$. 
Noting a cancellation of two terms,
and adding the same \gianni{quantity $\,\frac12\,\norma{\qn(t)}^2=\juerg{\frac 12\,\|q_n(T)\|^2\,}-\int_t^T\!\iO 
\qn\dt\qn\,$} to both sides, we arrive at the identity
\begin{align}
\label{disc4}
&\frac 12\,\|\pen(t)\|^2\,+\,\frac 12\,\|\qn(t)\|_{\VA\rho}^2\,+\int_t^T\!\!\!\iO|\dt\qn|^2\,+\int_t^T\!\!\!\iO|B^\sigma\pen|^2
\,+\pier{\int_t^T\!\!\!\iO F''_1(\bphi)\pen^2}
\nonumber\\
&=\,\frac 12\,\|\pen(T)\|^2\,+\,\frac 12\,\|\qn(T)\|^2_{\VA\rho}\,-\int_t^T\!\!\!\iO g_4\,\dt\qn \,+\int_t^T\!\!\!\iO g_2\pen
\,-\int_t^T\!\!\!\iO \pier{F''_2}(\bphi)\pen^2\nonumber\\
&\quad \,+\int_t^T\!\!\!\iO\ell'(\bphi)\btheta\,\pen^2-\int_t^T\!\!\!\iO\qn\dt\qn\,.
\end{align}
\pier{Note that the fifth term on the \lhs\ of \eqref{disc4} is nonnegative due to {\bf (F1)--(F2)}.}
By means of the H\"older and Young inequalities,
we readily conclude that the sum of the five integrals on the \rhs, which we denote by~$\,I$, satisfies
\begin{align}
\label{disc5}
|I|\,&\le\,\frac 12\int_t^T\!\!\!\iO|\dt\qn|^2\,+\,{C_1}\Bigl(\int_t^T\!\!\!\iO\bigl(|g_2|^2+|g_4|^2\bigr)\,+
\int_t^T\!\!\!\iO\bigl(|\qn|^2+|\pen|^2\bigr)\Bigr)
\nonumber\\
&\hspace*{5mm}+\,C_2\int_t^T\|\btheta(s)\|_{\Lx 4}\,\|\pen(s)\|\,\|\pen(s)\|_{\Lx 4}\,ds\nonumber\\
&\le \,\frac 12\int_t^T\!\!\!\iO|\dt\qn|^2\,+\,C_3\Bigl(\int_t^T\!\!\!\iO\bigl(|g_2|^2+|g_4|^2\bigr)\,+
\int_t^T\!\!\!\iO\bigl(|\qn|^2+|\pen|^2\bigr)\Bigr)\nonumber\\
&\hspace*{5mm}+\,\frac 12\int_t^T\!\!\!\iO \gianni{\bigl(|\pen|^2+|B^\sigma\pen|^2\bigr)}\,+\,
{C_4}\int_t^T\|\btheta(s)\|_{\VA\rho}^2\,\|\pen(s)\|^2\,ds\,.
\end{align}

It remains to estimate the final value terms. 
At first, note that the second identity in \eqref{disc3} just means that
$\pen(T)$ is the $H$-orthogonal projection of $g_1-\ell(\bphi(T))\,g_3$ onto~$V_n'$. 
Thus, $\,\|\pen(T)\|\le\|g_1\|+C_5\,\|g_3\|$. 
By the same token, \pier{we have that} $\|\qn(T)\|\le\|g_3\|$. 
Now observe that $\btheta(T)\in\VA\rho$. 
Therefore, invoking {\bf (A9)}, we have $g_3\in\VA\rho$. 
But this entails~that
\begin{align*}
\|A^\rho\qn \pier{(T)}\|^2&=(\qn\pier{(T)},A^{2\rho}\qn\pier{(T)})=(g_3,A^{2\rho}\qn\pier{(T)})\\
&=(A^\rho g_3,A^\rho\qn\pier{(T)}),\quad\mbox{i.e., }\,
\|A^\rho\qn\pier{(T)}\|\le\|A^\rho g_3\|.
\end{align*}
Hence, \pier{it turns out that} $\|\qn(T)\|_{\VA\rho}\le\|g_3\|_{\VA\rho}$. 
Observing that the mapping $\,\,s\mapsto\|\btheta(s)\|^2_{\VA\rho}\,\,$ belongs to $L^1(0,T)$, 
we obtain from the above estimates, using Gronwall's lemma, that
\begin{align}\label{disc6}
&\|\qn\|_{H^1(0,T;H)\cap L^\infty(0,T;\VA\rho)}\,+\,\|\pen\|_{L^\infty(0,T;H)\cap L^2(0,T;\VB\sigma)}\nonumber\\[0.5mm]
&\le\, C_6\left(\|g_1\|_{\Ldue}+\|g_2\|_{L^2(Q)}+\|g_3\|_{\VA\rho}+\|g_4\|_{L^2(Q)}\right)
\quad\forall\,n\in\enne.
\end{align}
Next, we insert $v=A^{2\rho}\qn$ in \eqref{disc1}. Using the estimate $\,\|\qn(T)\|_{\VA\rho}\le\|g_3\|_{\VA\rho}\,$ 
once more, we can infer that also 
\begin{equation}
\label{disc7}
\|\qn\|_{L^2(0,T;V_A^{2\rho})}\,\le\,C_7\left(\|g_1\|_{\Ldue}+\|g_2\|_{L^2(Q)}+\|g_3\|_{\VA\rho}+\|g_4\|_{L^2(Q)}\right) \quad\forall\,n\in\enne,
\end{equation}
and comparison in \eqref{disc2} shows that
\begin{equation}
\label{disc8}
\|\pen\|_{H^1(0,T;\VB{-\sigma})}\,\le\,C_8\left(\|g_1\|_{\Ldue}+\|g_2\|_{L^2(Q)}+\|g_3\|_{\VA\rho}+\|g_4\|_{L^2(Q)}\right)
\quad\forall\,n\in\enne.
\end{equation}
From the above estimates there  follows the existence of a pair $(q,p)$ such that, possibly only on a subsequence
which is still indexed by~$n$, 
\begin{align}\label{conq}
\qn\to q \quad&\mbox{weakly-star in }\,H^1(0,T;H)\cap L^\infty(0,T;\VA\rho)\cap  L^2(0,T;V_A^{2\rho}),\\[1mm]
\label{conp}
\pen\to p\quad&\mbox{weakly-star in }\,H^1(0,T;\VB {-\sigma})\cap L^\infty(0,T;H)\cap L^2(0,T;\VB\sigma). 
\end{align}
Moreover, by continuous embedding, $q\in C^0([0,T];\VA\rho)$ and $p\in C^0([0,T];H)$. 

At this point, it is a standard argument (which needs no repetition here) to show that $(q,p)$ is a 
solution to the system \eqref{as1}, \eqref{as2w}, \eqref{as3}. 
It remains to show uniqueness. 
To this end, let $(q_i,p_i)$, $i=1,2$, be two solutions, and $q=q_1-q_2$, $p=p_1-p_2$. 
Then $(q,p)$ solves \eqref{as1}, \eqref{as2w}, \eqref{as3} with zero right-hand sides. 
We now repeat the estimates leading to \eqref{disc6} for the continuous problem, concluding that $q=p=0$. 
The assertion is thus proved.
\Edim

We now can eliminate the variables $(\eta,\xi)$ from the variational inequality \eqref{vug1}.

\Bthm
\label{Optimality}
Let the assumptions {\bf (F1)}--{\bf (F4)}, {\bf (A1)}--{\bf (A9)}, and {\bf (GB)} be fulfilled, and let $\bu\in\uad$ be 
an optimal control of problem {\bf (CP)} with associated state $(\btheta,\bphi)={\cal S}(\bu)$ and adjoint state $(p,q)$. 
Then it holds the variational inequality 
\begin{equation}\label{vug2}
\iint_Q  q(u-\bu)\,+\beta_5\iint_Q\bu(u-\bu)\,\ge\,0\quad\forall\, u\in\uad\,. 
\end{equation}
\Ethm

\Bdim
\gianni{We fix $u\in\uad$ and consider the associated linearized system \eqref{ls1}--\eqref{ls3} with $h=u-\bu$}. 
We multiply
\eqref{ls1} by $\,q\,$ and \eqref{ls2} by $\,p\,$, add the results, and integrate over~$Q$.
\gianni{We obtain}
\begin{align}
\iint_Q q(u-\bu)\,&=\,
  \iint_Q \dt\eta \,q\,
  +\iint_Q \bigl( \ell'(\bphi)\,\dt\bphi\,\xi + \ell(\bphi)\,\dt\xi \bigr) \,q\, 
  +\iint_Q\eta\,A^{2\rho}q\,  
  \non
  \\
&\quad\,
  +\iint_Q \dt\xi \,p\,
  +\ioT\!\!(B^\sigma p(t),B^\sigma\xi(t))\,dt
\nonumber\\
\separa
&\quad\,
  +\iint_Q F''(\bphi)\,\xi\,p\,
  -\iint_Q\ell'(\bphi)\,\btheta\,\xi\,p\,
  - \iint_Q\ell(\bphi)\,\eta\,p \,.
\nonumber
\end{align}
\gianni{By also integrating by parts with respect to time in three of the terms, we deduce that}
\begin{align}
\iint_Q q(u-\bu)\,&=\,
  \iO \bigl(
    \eta(T)q(T)
    +\ell(\bphi(T))\xi(T)q(T)
    +\xi(T)p(T)
  \bigr)
\nonumber\\
\separa
&\quad\,
  +\iint_Q \eta\,\bigl[
    -\dt q+A^{2\rho}q
    -\ell(\bphi)p
  \bigr]\,
  -\ioT\!\!\langle\dt p(t),\xi(t)\rangle_{\VB\sigma}\,dt
\nonumber\\
\separa
&\quad\,
  +\ioT\!\!(B^\sigma p(t),B^\sigma\xi(t))\,dt\,
  +\iint_Q\xi\,\bigl[
    -\ell(\bphi)\dt q
    +F''(\bphi)p
    -\ell'(\bphi)\,\btheta\,p
  \bigr]\,.
  \non
\end{align}
Thus, using the adjoint system \eqref{as1}, \eqref{as2w}, \eqref{as3}, we find the identity
\begin{align*}
\iint_Q q(u-\bu)\,&=\beta_1\iO(\bphi(T)-\phi_\Omega)\xi(T)\,+\,\beta_2\iint_Q(\bphi-\phi_Q)\xi\\
&\quad\,+\,\beta_3\iO(\btheta(T)-\theta_\Omega)\eta(T)\,+\,\beta_4\iint_Q(\btheta-\theta_Q)\eta\,.
\end{align*}
\gianni{By combining this with \eqref{vug1}, we obtain~\eqref{vug2}.}
\Edim

\Brem
If $\beta_5>0$, then \eqref{vug2} just means that $\,\bu\,$ is the $L^2(Q)$-orthogonal projection of $\,-\beta_5^{-1}q\,$
onto $\uad$, i.e., we have 
\begin{equation}
\bu\,=\,\max\,\bigl\{u_{min},\,\min\,\bigl\{-\beta_5^{-1}q,u_{max}\bigr\}\bigr\} \quad\mbox{a.e.\ in }\,Q.
\end{equation}
\Erem


\subsection{The double obstacle case} 
\label{DQ}

In this section, we study the case of the double obstacle potential $F_{\rm 2obs}$ in which 
$F_1=I_{[-1,1]}$ is the indicator function of the interval $[-1,1]$ \pier{that} is given by 
\,$\indi(r)=0\,$ for $r\in [-1,1]$ \,and \,$\indi(r)=+\infty\,$ otherwise. 
Then the conditions
{\bf (F1)} and {\bf (F2)} are fulfilled with $(r_-,r_+)=(-1,1)$. 
For the other nonlinearities 
$F_2$ and $\ell$ we assume that {\bf (F3)} and {\bf (F4)} are valid. 
We then consider the following optimal control problem: 

\vspace{2mm}
\noindent
{\bf (CP$_0$)}  \,\,\,Minimize  \,\,${\cal J}((\theta,\phi),u)\,\,$ over \,$\uad$\, subject 
to the state system \eqref{ssvar1}--\eqref{ssvar3} with\linebreak
\hspace*{14mm} $\,F_1=\indi$.

\vspace{2mm}
\Brem
Notice that the condition $\,F_1(\phi)\in L^1(Q)\,$ for our notion of solution 
can \pier{only be satisfied for $\,F_1=\indi\,$} if $\,\phi\in[-1,1]$
almost everywhere, which in turn entails that the term involving
$F_1(\phi)$ on the \lhs\ of \eqref{ssvar2}  vanishes. 
\Erem

\vspace{2mm}
Since we cannot expect the condition {\bf (GB)} to be satisfied in this case, the control theory
developed in the previous section does not apply. 
We therefore argue by approximation,
using the deep quench approximation, which has proved to be successful in a number of 
similar situations (see, e.g., \cite{CFS,CFGS,CGS_DQ,CGS21, CGS21bis,S_DQ}). 
The general idea behind this approach is the following: we define the logarithmic functions
\begin{align}
\label{defh}
&h(r):=\left\{
\begin{array}{l}
(1+r)\ln(1+r)+(1-r)\ln(1-r)\quad\mbox{if $\,r\in (-1,1)$}\\[0.5mm]
2\ln(2)\quad\mbox{if \,$r\in\{-1,1\}$}\\[0.5mm]
+\infty\quad\mbox{if \,$r\not\in [-1,1]$}
\end{array}
\right.\\[1mm]
\label{defhal}
&h_\alpha(r):=\alpha h(r)\quad\mbox{for \,$r \in\erre$\, and \,$\alpha\in (0,1]$}.
\end{align}
It is easily seen that
\Beq
\label{halcon}
\lim_{\alpha\searrow0}\,\hal(r)=\indi(r)\quad\forall\,r\in\erre .
\Eeq
Moreover, $h'(r)=\ln(\frac{1+r}{1-r})$\, and \,$h''(r)=\frac 2{1-r^2}$, and thus
\begin{align}
\label{halscon}
&\lim_{\alpha\searrow0} \dhal(r)=0 \quad\,\hbox{for all }\,r\in (-1,1),
\nonumber\\
&\quad \lim_{\alpha\searrow0}\Bigl(\lim_{r\searrow -1}\,\dhal (r)\Bigr)=-\infty, \quad
\lim_{\alpha\searrow0}\Bigl(\lim_{r\nearrow 1}\,\dhal (r)\Bigr)=+\infty.
\end{align}
Hence, we may regard the graphs of the single-valued functions $\dhal$ over 
the interval $(-1,1)$ as approximations to the graph of the subdifferential~$\partial\indi$.
Observe that this is an \emph{interior} approximation defined in the interior of the 
domain of $\partial\indi$ in contrast to the \emph{exterior} approximation obtained via
the Moreau--Yosida approach. 

In view of \eqref{halcon}--\eqref{halscon}, it is near to mind to expect that the control problem
{\bf (CP$_0$)} is closely related to the control problem
(which in the following will be denoted by {\bf (CP$_\alpha$)}) 
that arises when in  \eqref{ssvar2} we choose $F_1=\hal$ for $\alpha>0$. 
Indeed, by virtue of Theorem~\ref{FromCG}, the system \eqref{ssvar1}--\eqref{ssvar3} 
enjoys for both $F_1=\indi$ and $F_1=\hal$ a solution pair $(\theta,\phi)$ and $(\tal,\pal)$.
We introduce the corresponding solution operators 
$${\cal S}_0:{\cal U}_R\ni u\mapsto (\theta,\phi),\quad {\cal S}_{\alpha}:{\cal U}_R
\ni u\mapsto (\tal,\pal).$$ 
It can be expected that  $(\tal,\pal)$ converges in a suitable
topology to $(\theta,\phi)$ as $\alpha\searrow0$. 
Moreover,  the optimal
control problem {\bf (CP$_\alpha$)} belongs to the class of problems for which in Section~\ref{NECOPT} 
first-order necessary optimality conditions in terms of a variational inequality 
and the adjoint state system have been established. 
One can therefore hope to perform a passage to the limit as $\alpha\searrow 0$ in the state and the adjoint state
variables in order to derive meaningful first-order necessary optimality conditions also for~{\bf (CP$_0$)}.

In order to carry out this program, we now make a restrictive assumption, which still includes
the classical situation: 

\vspace{2mm} \noindent
\gianni{{\bf (A10)} \,\,\,It holds $\,B^{2\sigma}=B=-\Delta\,$ with zero Dirichlet or Neumann boundary conditions,\linebreak 
\hspace*{14mm} and $A=-\Delta$ with
zero Neumann or Dirichlet boundary conditions with either \linebreak \hspace*{14mm} $\rho>\frac34$ or $\rho=1/2$.}

\vspace{2mm}
\Brem
\label{FromA10}
If {\bf (A10)} is valid, then the conditions {\bf (A4)}, {\bf (A7)} and {\bf (A8)} are 
automatically satisfied.   
\Erem

We now assume that also the assumptions {\bf (A1)}--{\bf (A3)}, {\bf (A5)} and {\bf (A6)} are fulfilled. 
Now observe that under {\bf (A10)} 
\gianni{both the assumption~(i) of Lemma~\ref{SuffGB} and the condition \eqref{partin} are met 
(see Remark~\ref{Rempartin})}. 
\gianni{Since the functions $\hal$ satisfy the condition \eqref{asymp}, we thus can conclude
from Lemma~\ref{SuffGB}(i) and its proof that 
the solutions $(\tal,\pal)$ to the state system with $F_1=\hal$ satisfy both the boundedness condition \eqref{yeah}
and the condition {\bf (GB)} for every $\alpha>0$. 
Therefore}, for every $\alpha>0$, there are constants $a^\alpha_R,b^\alpha_R\gianni{,c^\alpha_R}$ such that
\begin{equation}\label{separal}
-1<a^\alpha_R\le\pal\le b^\alpha_R<1
\gianni{\aand
|\tal|\leq c^\alpha_R}
\quad\mbox{a.e.\ in }\,Q,
\end{equation}
whenever $(\tal,\pal)={\cal S}_{\alpha}(u)$ for some $u\in{\cal U}_R$. 
\pier{In addition,} as it was \gianni{established in Lemma~\ref{FromA7}}, the variational inequality \eqref{ssvar2}
takes for every $\alpha>0$ the form of a variational equality, namely
\gianni{%
\begin{align}
\label{ssvar2al}
&(\dt\phi_\alpha(t),v)+(\nabla\phi_\alpha(t),\nabla v)+(\dhal(\phi_\alpha(t)),v)+(F_2'(\phi_\alpha(t)),v)
\nonumber\\
&\,=\,(\ell(\phi_\alpha(t))\theta_\alpha(t),v)
\quad\mbox{for a.e.\ $t\in(0,T)$ and all }\,v\in H^1(\Omega),
\end{align}
and $(\tal,\pal)$ is in fact a strong solution}. 

The approximating control problem reads:

\vspace{2mm}\noindent
{\bf (CP$_\alpha$)} \,\,\,Minimize the cost functional \eqref{cost} over $\uad$ subject 
to the state system
\eqref{ssvar1},\linebreak\hspace*{15mm} \eqref{ssvar2al}, \eqref{ssvar3}.

\vspace{3mm}
We \gianni{recall Remark~\ref{FromA10} and state} the following approximation result.

\Bthm
\label{Boundal}
Suppose that {\bf (F3)}, {\bf (F4)}, {\bf (A1)}--{\bf (A3)}, {\bf (A5)}--{\bf (A6)},
and {\bf (A10)} are fulfilled, and assume that $(\tal,\pal)={\cal S}_\alpha(u_\alpha)$ 
for some $u_\alpha\in{\cal U}_R$ and $\alpha\in (0,1]$. 
Then there is some constant $K_3>0$, which depends only on $R$ and the data, such that 
\begin{align}
\label{ssbound1al}
&\|\tal\|_{H^1(0,T;H)\cap L^\infty(0,T;\VA\rho)\cap L^2(0,T;\VA{2\rho})}\nonumber\\
&\quad +\,
\|\pal\|_{W^{1,\infty}(0,T;H)\cap H^1(0,T;\Huno)}+\iint_Q \hal \pier{(\pal)}\,\le\,K_3.
\end{align} 
Moreover, there is a sequence $\{\alpha_n\}\subset (0,1]$ with $\alpha_n\searrow0$ such that  
\begin{align}\label{convual}
u_{\alpha_n}\to u\,\quad&\mbox{weakly-star in }\,L^\infty(Q),\\[1mm]
\label{convtal}
\taln\to\theta\,\quad&\mbox{weakly-star  in }\, H^1(0,T;H)\cap L^\infty(0,T;\VA\rho)\cap 
L^2(0,T;\VA{2\rho}),\nonumber\\[0.5mm]
&\mbox{strongly in \,$C^0([0,T];H)\,$ and pointwise a.e.\ in $\,Q$,}\\[1mm]
\label{convpal}
\paln\to\phi\,\quad&\mbox{weakly-star in }\,W^{1,\infty}(0,T;H)\cap H^1(0,T;\Huno),\nonumber\\[0.5mm]
&\mbox{strongly in \,$C^0([0,T];H)\,$ and pointwise a.e.\ in $\,Q$,}
\end{align}
where $(\theta,\phi)$ denotes the unique solution to the state system 
\eqref{ssvar1}--\eqref{ssvar3} for $F_1=\indi$ and the control~$\,u$.
\Ethm

\Bdim
The validity of the estimate \eqref{ssbound1al} follows from a closer inspection of the
derivation of the a priori estimates performed in~\cite{CG}: 
indeed, \pier{by virtue of {\bf (A3)}} it turns out that
the bounds derived there are for $F_1=\hal$ in fact independent of $\alpha\in (0,1]$. 
Hence, there are a sequence $\alpha_n\searrow0$ and  $u,\theta,\phi$ satisfying 
\eqref{convual}--\eqref{convpal}, where the strong convergence in $C^0([0,T];H)$ follows from
\cite[Sect.~8, Cor.~4]{Simon}. 
It remains to show that $(\theta,\phi)={\cal S}_0(u)$. 

At first, it is easily seen that $\theta(0)=\theta_0$ and $\phi(0)=\phi_0$. 
Moreover, we observe that \eqref{convpal} entails, by Lipschitz continuity, that 
$\,\ell(\paln)\to\ell(\phi)\,$ and $\,F_2'(\paln)\to F_2'(\phi)$, both strongly in $C^0([0,T];H)$. 
Moreover, the sequences $\{\ell(\paln)\dt\paln\}$ and $\{\ell(\paln)\taln\}$ are bounded in $L^2(Q)$
\gianni{since $\ell$ is bounded}. 
This entails that
$$
\ell(\paln)\dt\paln\to\ell(\phi)\dt\phi\quad\mbox{and}\quad\ell(\paln)\taln\to
\ell(\phi)\theta, \quad\mbox{both weakly in }\,L^2(Q).
$$ 
Hence, we may write \eqref{ssvar1}, with $F_1=\haln$  and control $\ualn$, and pass to the
limit as $n\to\infty$ to see that $(\theta,\phi)$ satisfies \eqref{ssvar1} with control~$\,u$.
It remains to show \eqref{ssvar2} with $F_1=\indi$. 
We are going to prove it in the time-integrated form \eqref{intssvar2}. 

To this end, we first note that \eqref{ssbound1al} entails that we must have $\,\paln\in[-1,1]\,$ a.e.\ in $\,Q\,$. 
Since $\,\paln\to\phi\,$ pointwise a.e.\ in $\,Q$, also $\phi\in[-1,1]$ a.e.\ in $\,Q\,$ and thus $\iint_Q\indi(\phi)=0$.

Now let $v\in L^2(0,T;\Huno)$ be arbitrary. 
If $\,\indi(v)\not\in L^1(Q)$, then the inequality is
fulfilled since its \rhs\ is infinite. 
Otherwise, we have $\,v\in[-1,1]\,$ a.e.\ in $\,Q\,$ and thus 
 $0=\indi(v)\le\haln(v)\le h_1(v)$ a.e.\ in $\,Q$. Since, thanks to \eqref{halcon},
 $\,\haln(v)\to \indi(v)\,$ pointwise a.e.\ in $\,Q$, we infer from Lebesgue's dominated convergence
 theorem that $\,0=\iint_Q\indi(v)=\lim_{n\to\infty} \iint_Q\haln(v)$.
Therefore, using the lower semicontinuity of the quadratic \gianni{form $\,v\mapsto\iint_Q|\nabla v|^2\,$}
on $L^2(0,T;\Huno)$, we can infer that
\begin{align*}
&\iint_Q\indi(\phi)+\iint_Q \nabla\phi\cdot\nabla(\phi-v)\,\le\,\liminf_{n\to\infty}\iint_Q\nabla\paln\cdot\nabla(\paln-v)
\\
\separa
&\le\,\liminf_{n\to\infty}\,\Bigl(\iint_Q\bigl(\ell(\paln)\taln-\dt\paln-F_2'(\paln)\bigr)(\paln-v)\,+\iint_Q\haln(v)\Bigr)
\\
&\gianni{{}={}}\iint_Q\bigl(\ell(\phi)\theta-\dt\phi-F_2'(\phi)\bigr)(\phi-v)\,+\iint_Q \indi(v)\,.
\end{align*}
This finishes the proof of the assertion.
\Edim

\Brem
Notice that a uniform (with respect to $\,\alpha\in (0,1]$) bound resembling \eqref{ssbound2} for $F_1=\hal$ cannot be expected to hold true, 
since it may well happen that $\,a^\alpha_R\searrow -1\,$ and/or $\,b^\alpha_R\nearrow +1\,$ as $\alpha\searrow0$, 
so that $\dhal(\pal)$ and $h_\alpha''(\pal)$ may become unbounded as $\alpha\searrow0$.
\Erem

\vspace{2mm}
\pier{In view of the expression \eqref{cost} of the functional $\cal J$ and of Theorem~\ref{Boundal},
it is not difficult to argue that}
 optimal controls of {\bf (CP$_\alpha$)} 
are ``close'' to optimal controls of {\bf (CP$_0$)}. 
However, \pier{from Theorems~\ref{Boundal} we cannot infer sufficient information 
on the family of} the minimizers of~{\bf (CP$_0$)}. 
In order to find first-order necessary optimality conditions, 
we recall that in the previous section we have been able to derive such conditions for the problem~{\bf (CP$_\alpha$)}. 
Thus, we can hope to establish corresponding results for {\bf (CP$_0$)} by taking the limit as $\alpha\searrow0$. 
However, such an approach fails since the convergence property
\eqref{convual} is too weak to pass to the limit as $\alpha\searrow0$ in the variational inequality
\eqref{vug2} (written for an optimal control $\overline u_\alpha$ and the corresponding adjoint state~$q_\alpha$). 
For this, we seem to need a strong convergence of $\{\overline u_\alpha\}$ in~$L^2(Q)$.  

To this end, we employ a well-known technique. 
Let us assume that
$\bu\in\uad$ is any optimal control for {\bf (CP$_0$)} with associated state $(\btheta,\bphi)={\cal S}_0(\bu)$. 
We associate with it the {\em adapted cost functional}
\begin{equation}
\label{adcost}
\widetilde{\cal J}((\theta,\phi),u):={\cal J}((\theta,\phi),u)\,+\,\frac 12\,\|u-\bu\|^2_{L^2(Q)}
\end{equation}
and a corresponding \emph{adapted optimal control problem}:

\vspace{2mm}\noindent
{\bf  ($\widetilde{\mathbf{CP}}_\alpha$)} \,\,\,Minimize the cost functional \eqref{adcost} over $\uad$ subject to the state system
\eqref{ssvar1}--\linebreak \hspace*{16mm}\eqref{ssvar3}, where $F_1=h_\alpha$.

\vspace{2mm}
With the same direct argument as in the proof of Theorem~\ref{Existoc}, 
we can show that {\bf  ($\widetilde{\mathbf{CP}}_\alpha$)} has a solution. 
The following result indicates why the 
adapted control problem suits better for our intended approximation approach.

\Bthm
\label{Convadapt}
Suppose that {\bf (F3)}\pier{--}{\bf (F4)}, {\bf (A1)}--{\bf (A3)}, {\bf (A5)}--{\bf (A6)},
and {\bf (A10)} are fulfilled, assume that 
$\bu\in \uad$ is an arbitrary optimal control of {\bf (CP$_0$)} with associated state  
$(\btheta,\bphi)$, and let $\,\{\alpha_n\}\subset (0,1]\,$ be
any sequence such that $\,\alpha_n\searrow 0\,$ as $\,n\to\infty$. 
Then there exist a subsequence $\{\alpha_{n_k}\}_{k\in\enne}$ of $\{\alpha_n\}$, and, for every $k\in\enne$, 
an optimal control $\,\bu_{\alpha_{n_k}}\in \uad\,$ of the adapted problem {\bf ($\widetilde{\mathbf{CP}}_{\alpha_{n_k}}$)}
with associated state $(\btheta_{\alpha_{n_k}},\bphi_{\alpha_{n_k}})$
such that, as $k\to\infty$,
\begin{align}
\label{adap1}
&\bu_{\alpha_{n_k}}\to \bu\quad\mbox{strongly in }\,L^2(Q),
\end{align}
and the properties  \eqref{convtal} and \eqref{convpal} are satisfied correspondingly. Moreover, we have 
\begin{align}
\label{adap2}
&\lim_{k\to\infty}\,\widetilde{{\cal J}}((\btheta_{\alpha_{n_k}},\bphi_{\alpha_{n_k}}),
\bu_{\alpha_{n_k}})\,=\,  {\cal J}((\btheta,\bphi),\bu)\,.
\end{align}
\Ethm

\Bdim
Let $\alpha_n \searrow 0$ as $n\to\infty$. 
For any $ n\in\enne$, we pick an optimal control 
$u_{\alpha_n} \in \uad\,$ for the adapted control problem {\bf ($\widetilde{\mathbf{CP}}_{\alpha_{n}}$)}
and denote by 
$(\theta_{\alpha_n},\phi_{\alpha_n})$ the associated solution to  the state system 
with $F_1=h_{\alpha_n}$ and $u=u_{\alpha_n}$. 
By the boundedness of $\uad$ in $L^\infty(Q)$, there is some subsequence $\{\alpha_{n_k}\}$ of $\{\alpha_n\}$ such that
\begin{equation}
\label{ugam}
u_{\alpha_{n_k}}\to u\quad\mbox{weakly-star in }\,L^\infty(Q)
\quad\mbox{as }\,k\to\infty,
\end{equation}
with some $u\in\uad$\pier{, and, thanks} to Theorem~\ref{Boundal}, 
the convergence properties \eqref{convtal} and \eqref{convpal} hold true 
with the pair $(\theta,\phi)={\cal S}_0(u)$. 
In particular,  the pair $((\theta,\phi),u)$ is admissible for~{\bf (CP$_0$)}. 

We now aim to prove that $u=\bu$. 
Once this is shown, it follows from the unique solvability of the state system that also $(\theta,\phi)=(\btheta,\bphi)$, 
which implies that \eqref{convtal} and \eqref{convpal} hold true with  $\,(\theta,\phi)\,$ replaced by $(\btheta,\bphi)$.

Now observe that, owing to the weak sequential lower semicontinuity of~$\widetilde {\cal J}$, 
and in view of the optimality property of $((\btheta,\bphi),\bu)$ for problem~{\bf (CP$_0$)},
\begin{align}
\label{tr3.6}
&\liminf_{k\to\infty}\, \widetilde{\cal J}((\theta_{\alpha_{n_k}},\phi_{\alpha_{n_k}}),
u_{\alpha_{n_k}})
\ge \,{\cal J}((\theta,\phi),u)\,+\,\frac{1}{2}\,
\|u-\bu\|^2_{L^2(Q)}\nonumber\\[1mm]
&\geq \, {\cal J}((\btheta,\bphi),\bu)\,+\,\frac{1}{2}\,\|u-\bu\|^2_{L^2(Q)}\,.
\end{align}
On the other hand, the optimality property of  $\,((\theta_{\alpha_{n_k}},\phi_{\alpha_{n_k}}),u_{\alpha_{n_k}})\,$ 
for problem {\bf ($\widetilde{\mathbf{CP}}_{\alpha_{n_k}}$)} yields that for any $k\in\enne$ we have
\begin{equation}
\label{tr3.7}
\widetilde {\cal J}((\theta_{\alpha_{n_k}},\phi_{\alpha_{n_k}}),u_{\alpha_{n_k}})\, =\,
\widetilde {\cal J}({\cal S}_{\alpha_{n_k}}(u_{\alpha_{n_k}}),u_{\alpha_{n_k}})
\,\le\,\widetilde {\cal J}({\cal S}_{\alpha_{n_k}}(\bu),\bu)\,.
\end{equation}
\gianni{Finally, with the same argument used at the beginning of the proof of Theorem~\ref{Boundal}
to justify the estimate~\eqref{ssbound1al}, 
one sees that ${\cal S}_{\alpha_{n_k}}(\bu)$ satisfies a similar bound,
whence \pier{a subsequence (not relabeled)} converges to some pair~$(\underline\theta,\underline\phi)$ in the topologies specified in \accorpa{convtal}{convpal}.
As in the proof of the abovementioned theorem, one shows that $(\underline\theta,\underline\phi)$ 
solves the original state system associated with~$\bu$, i.e., it coincides with~${\cal S}_0(\bu)$.
Therefore, invoking the continuity properties
of the cost functional with respect to the topologies of the spaces $C^0([0,T];H)$ and $L^2(Q)$, we deduce from \eqref{tr3.7} that}
\gianni{%
\begin{align}
\label{tr3.8}
&\limsup_{k\to\infty}\,\widetilde {\cal J}((\theta_{\alpha_{n_k}},\phi_{\alpha_{n_k}}),u_{\alpha_{n_k}})
\,\le\,\limsup_{k\to\infty}\,\widetilde {\cal J}({\cal S}_{\alpha_{n_k}}(\bu),\bu)
\non
\\
&\,=\,\widetilde {\cal J}({\calS}_0(\bu),\bu) 
\,=\,\widetilde {\cal J}((\btheta,\bphi),\bu)
\,=\,{\cal J}((\btheta,\bphi),\bu)\,.
\end{align}
}%
Combining (\ref{tr3.6}) with (\ref{tr3.8}), we have thus shown that 
$\,\frac{1}{2}\,\|u-\bu\|^2_{L^2(Q)}=0$\,,
so that $\,u=\bu\,$  and thus also $\,(\theta,\phi)=(\btheta,\bphi)$. 
Moreover, (\ref{tr3.6}) and (\ref{tr3.8}) also imply that
\begin{align*}
&{\cal J}((\btheta,\bphi),\bu) \, =\,\widetilde{\cal J}((\btheta,\bphi),\bu)
\,=\,\liminf_{k\to\infty}\, \widetilde{\cal J}((\theta_{\alpha_{n_k}},\phi_{\alpha_{n_k}}),
 u_{\alpha_{n_k}})\nonumber\\[1mm]
&\,=\,\limsup_{k\to\infty}\,\widetilde{\cal J}((\theta_{\alpha_{n_k}},\phi_{\alpha_{n_k}}),
u_{\alpha_{n_k}}) \,
=\,\lim_{k\to\infty}\, \widetilde{\cal J}((\theta_{\alpha_{n_k}},\phi_{\alpha_{n_k}}),
u_{\alpha_{n_k}})\,,
\end{align*}                                     
which proves \pier{\eqref{adap1} and {(\ref{adap2})} at the same time, of course along with} 
(\ref{convtal}) and \eqref{convpal}. 
This concludes the proof of the assertion.
\Edim

\vspace{2mm}
We now discuss the first-order necessary optimality conditions for~{\bf ($\widetilde{\mathbf{CP}}_{\alpha}$)}, 
assuming that the general assumptions {\bf (F3)}\pier{--}{\bf (F4)}, {\bf (A1)}--{\bf (A3)}, {\bf (A5)}--{\bf (A6)},
{\bf (A9)} and {\bf (A10)} are fulfilled. 
Obviously, the adjoint system is the same as for~{\bf (CP$_\alpha$)}, 
and Theorem~\ref{Wellposadj} and Theorem~\ref{Optimality} apply to this situation.
More precisely, the adjoint state $(p_\alpha,q_\alpha)$ solves the variational system 
\begin{align}
\label{as1al}
&-\dt q_\alpha-\ell(\bphi_\alpha)\,p_\alpha+A^{2\rho}q_\alpha\,=\,g_4^\alpha\quad\mbox{in }\,Q,\\[0.5mm]
\label{as2al}
&\gianni{\bigl( -\dt p_\alpha(t),v \bigr)}
  \,-\,\bigl( \ell(\bphi_\alpha(t))\,\dt q_\alpha(t),v \bigr)
  +(\nabla p_\alpha(t),\nabla v)
\nonumber\\
&\quad {}
  +((\psi_1^\alpha(t)+\psi_2^\alpha(t))\,p_\alpha(t),v)
  -(\ell'(\bphi_\alpha(t))\btheta_\alpha(t)p_\alpha(t),v)\,=\,(g_2^\alpha(t),v)
\non\\
&\quad \, \mbox{for all $\,v\in\Huno\,$ and a.e.\ $t\in(0,T)$},\\[0.5mm]
\label{as3al}
&q_{\alpha}(T)=\pier{g_3^\alpha} , \quad p_\alpha(T)=g_1^\alpha-
\ell(\bphi_\alpha(T))\,\pier{g_3^\alpha}\quad
\mbox{in }\,\Omega,
\end{align}
where, for $\alpha>0$,
\begin{align}\label{notal}
  &\psi_1^\alpha := h_\alpha''(\bphi_\alpha), \quad \psi_2^\alpha :=F_2''(\bphi_\alpha) , 
  \quad g_1^\alpha := \beta_1(\bphi_\alpha(T)-\phi_\Omega), \quad
    g_2^\alpha := \beta_2(\bphi_\alpha-\phi_Q),\nonumber\\
  &g_3^\alpha:=\beta_3(\btheta_\alpha(T)-\theta_\Omega),
  \quad g_4^\alpha:=\beta_4\gianni{(\btheta_\alpha-\theta_Q)}\,.
\end{align}
By virtue of  the general bounds \eqref{separal}, \pier{\eqref{ssbound1al}} and owing to~{\bf (A9)}, we have that
\Beq
\label{restal}
\|\psi_2^\alpha\|_{L^\infty(Q)}+\pier{\|g_1^\alpha\|}
+\|g_2^\alpha\|_{L^2(Q)}\,+\|g_3^\alpha\|_{\VA\rho}+\|g_4^\alpha\|_{L^2(Q)}\,\le\, C_1 \quad\forall \,\alpha\in (0,1],
\Eeq
where, here and in the following, $C_i>0$, $i\in\enne$, denote constants that may depend on the data of 
the system, but not on $\alpha\in (0,1]$. 
Observe that a corresponding bound for $\psi^\alpha_1$ cannot be expected.

On the other hand, the variational inequality characterizing optimal controls 
is different (nevertheless, obtained using the same arguments that led to \eqref{vug2} in Theorem~\ref{Optimality}). 
Namely, 
if $\bu_\alpha \in\uad$ is optimal for {\bf  ($\widetilde{\mathbf{CP}}_\alpha$)} 
and $(p_\alpha,q_\alpha)$ is the associated
adjoint state, then we have that
\begin{align}
\label{vugal}
 \iint_Q (q_\alpha+\beta_5\bu_\alpha+(\bu_\alpha-\bu))(u-\bu_\alpha)\,\ge\,0 
  \quad\forall\,u\in\uad.
\end{align}
Our aim is to let $\alpha$ tend to zero in both the above inequality and the adjoint system.
Thus, we have to derive some a priori estimates for the adjoint variables that are uniform with respect to $\alpha\in (0,1]$. 
To this end, we note that the estimates \eqref{disc6}, \eqref{disc7}, derived for the Faedo--Galerkin approximations,
persist by the semicontinuity of norms under limit processes, whence we infer that
\begin{align}
\label{pqboundal}
&\|q_\alpha\|_{H^1(0,T;H)\cap L^\infty(0,T;\VA\rho)\cap L^2(0,T;\VA{2\rho})}\,+\,
\|p_\alpha\|_{L^\infty(0,T;H)\cap L^2(0,T;\Huno)}\nonumber\\[0.5mm]
&\le \,C_2\left(\pier{\|g_1^\alpha\|}
+\|g_2^\alpha\|_{L^2(Q)}\,+\|g_3^\alpha\|_{\VA\rho}+\|g_4^\alpha\|_{L^2(Q)}\right)\,\le\,C_3\quad\forall\,\alpha\in (0,1].
\end{align} 
However, the comparison argument leading to \eqref{disc8} does not work in this situation, 
because we do not have a bound for~$\psi_1^\alpha$. 
For this reason, we introduce the space
\begin{align}
\label{defZ}
{\cal Z}:=\{v\in H^1(0,T;\Huno^*)\cap L^2(0,T;\Huno):\,v(0)=0\}.
\end{align}
Since the embedding $(H^1(0,T;\Huno^*)\cap L^2(0,T;\Huno))\subset C^0([0,T];H)$ is continuous, the zero condition
for the initial value is meaningful, and $\,{\cal Z}\,$ is a closed subspace of $H^1(0,T;\gianni{\Huno^*})\cap L^2(0,T;\Huno)$
and thus a Banach space when endowed with the natural norm of this space. 
Moreover, the embedding $\,{\cal Z}\subset C^0([0,T];H)\,$ is continuous, and we also have the dense and 
continuous embedding ${\cal Z}\subset L^2(0,T;H)\subset{\cal Z}^*$,
where it is understood that 
\begin{equation}
\label{dualZ}
\langle v,z\rangle_{\cal Z}=\int_0^T(v(t),z(t))\,dt \quad\mbox{for all \gianni{\,$v\in L^2(0,T;H)$ and $\,z\in {\cal Z}\,$.}}
\end{equation}

\pier{Now, let $v\in {\cal Z}$ be arbitrary and use it as test function in \eqref{as2al}.  By integrating  \eqref{as2al} over $(0,T)$, with the help of~\eqref{as3al} we find out that 
\begin{align}
&\iO  \bigl( \ell(\bphi_\alpha(T))\,\pier{g_3^\alpha} - g_1^\alpha \bigr) v(T) + \ioT \langle   \dt v (t) , p_\alpha (t) \rangle_{H^1(\Omega)} dt
\nonumber\\
&{}-\iint_Q  \ell(\bphi_\alpha)\,\dt q_\alpha v + \iint_Q \nabla p_\alpha \cdot \nabla v + \iint_Q 
\psi_1^\alpha \,p_\alpha \, v   
\nonumber \\
&{}+ \iint_Q \psi_2^\alpha \,p_\alpha v - \iint_Q \ell'(\bphi_\alpha)\,\btheta_\alpha\, p_\alpha\, v 
=  \iint_Q g_2^\alpha\, v  . 
\label{pier1}
\end{align}
Then,  due to {\bf (F4)}, \eqref{restal} and \eqref{pqboundal}, we have that 
\begin{align}
\label{1a1}
&\Bigl| \iO  \bigl( \ell(\bphi_\alpha(T))\,\pier{g_3^\alpha} - g_1^\alpha \bigr) v(T)+ \ioT \langle   \dt v (t) , p_\alpha (t) \rangle_{H^1(\Omega)} dt \Bigr| \nonumber\\
&\le\,  \bigl(  C_3 \juergen{\|g_3^\alpha\|\,+\, \|g_1^\alpha\|}\bigr) \norma v_{\C0H}  \nonumber\\
&\quad\, +  \norma{p_\alpha}_{L^2(0,T;H^1(\Omega))}\,\norma{\dt v}_{L^2(0,T;H^1(\Omega)^* )}\, \le\,C_4\,\|v\|_{{\cal Z}}\,.
\end{align}
Moreover, \eqref{pqboundal} obviously yields that
\begin{align}\label{2al}
\| \pier{{} \ell(\bphi_\alpha)\,\dt q_\alpha}\|_{L^2(Q)} + \| \psi_2^\alpha\,\gianni{p_\alpha}\|_{L^2(Q)} + \| g_2^\alpha\|_{L^2(Q)}
\,\le\,\gianni{C_5} . 
\end{align}}%
\gianni{Furthermore, using also \eqref{ssbound1al} for $(\btheta_\alpha,\bphi_\alpha)$ and the fact that {\bf (A10)} 
implies~{\bf(A8)} (see Remark~\ref{FromA10}) \juerg{and thus also the continuity of the} embedding $\VA{2\rho}\subset\Lx4$, we \pier{deduce}~that}
\begin{align}\label{3al}
&\Bigl|\iint_Q \nabla\gianni{p_\alpha}\cdot\nabla v
\,-\iint_Q\ell'(\bphi_\alpha)\btheta_\alpha\,\gianni{p_\alpha}\,v\Bigr|\nonumber\\
&\le\,\gianni{C_6}\,\|v\|_{L^2(0,T;\Huno)}\,+\,\gianni{C_7}\ioT\|\btheta_\alpha(t)\|_{\Lx 4}\,\|\gianni{p_\alpha}(t)\|_{\Lx 4}\,\|v(t)\|\,dt\nonumber\\[1mm]
&\le\,\gianni{C_6}\,\|v\|_{\cal Z}\,+\,\gianni{C_8}\,\|\btheta_\alpha\|_{L^2(0,T;\VA{2\rho})}\,\|\gianni{p_\alpha}\|_{L^2(0,T;\Huno)}\,\pier{\|v\|_{C^0([0,T];H)}}
\,\le\,C_9\,\|v\|_{\cal Z}\,
\end{align}
for all $v\in{\cal Z}$. 
Hence, comparison in \pier{\eqref{pier1}} leads to 
\gianni{%
\begin{equation}
\label{lamal}
\|\Lambda_\alpha\|_{{\cal Z}^*}\,\le\,C_{10} \,,
\quad \hbox{with $\,\Lambda_\alpha:=\psi_1^\alpha p_\alpha=\alpha\,h''(\bphi_\alpha) p_\alpha\,$},
 \quad\forall\,\alpha\in (0,1].
\end{equation}}%

At this point, we are in a position to show the following first-order optimality result.

\Bthm
\label{NC2obs}
Suppose that the conditions {\bf (F3)}\pier{--}{\bf (F4)}, {\bf (A1)}--{\bf (A3)}, {\bf (A5)}--{\bf (A6)},
{\bf (A9)} and {\bf (A10)} are fulfilled, and let $\,\bu\in\uad\,$ be an optimal control
for {\bf (CP$_0$)} with associated state $\,(\btheta,\bphi)$. 
Then there exist $(q,p,\Lambda)$ such that the following statements hold true:
\gianni{%
\noindent\hbox to 9mm{{\rm (i)}\hfil}$q\in H^1(0,T;H)\cap C^0([0,T];\VA\rho)\cap L^2(0,T;\VA{2\rho})$,
\\[0.5mm]
\noindent\hbox to 9mm{\hfil}$p\in \pier{\L\infty H \cap L^2(0,T;\Huno)}$,%
\aand
$\Lambda \in{\cal Z}^*$.
}%

\gianni{%
\noindent\hbox to 9mm{{\rm (ii)}\hfil}The adjoint system, consisting of \eqref{as1}, \pier{the final condition 
\begin{equation}
q(T)=\beta_3(\btheta(T)-\theta_\Omega) \quad
\mbox{in }\,\Omega
\label{pier2}
\end{equation}
\noindent\hbox to 9mm{\hfil}and the equation
\begin{align}
&\iO  \bigl(\beta_3 \ell(\bphi(T))(\btheta(T)-\theta_\Omega)  - \beta_1 (\bphi(T)-\phi_\Omega)  \bigr) v(T) 
\nonumber\\
&{}+ \ioT \langle   \dt v (t) , p(t) \rangle_{H^1(\Omega)} dt-\iint_Q \ell(\bphi)\,\dt q\,v
+\iint_Q\nabla p\cdot \nabla v
+\langle \Lambda,v\rangle_{{\cal Z}} 
\nonumber \\
&{}+\iint_Q F_2''(\bphi)\,p\,v -\iint_Q\ell'(\bphi))\btheta\,p\,v
\,=\,\beta_2\iint_{\juerg{Q}}(\bphi-\phi_Q)\,v\quad \hbox{for all } v \in \cal Z,  \label{limas2}
\end{align}
\noindent\hbox to 9mm{\hfil}is satisfied.}%
}%

\gianni{%
\noindent\hbox to 9mm{{\rm (iii)}\hfil}It holds the variational inequality                           
\begin{align}
\label{vug3}
\iint_Q(q\,+\,\beta_5\,\bu)(u-\bu)\,\ge 0\quad \hbox{for all } u\in\uad \,.
\end{align}
}%
\Ethm

\Bdim
We choose any sequence $\{\alpha_n\}$ such that $\alpha_n\searrow0$. 
By Theorem~\ref{Convadapt} we may assume that there are
optimal controls  $\,\bu_{\alpha_n}\in \uad\,$ of the adapted problem {\bf ($\widetilde{\mathbf{CP}}_{\alpha_{n}}$)}
with associated states $(\btheta_{\alpha_{n}},\bphi_{\alpha_{n}})$
such that \eqref{adap1} \gianni{and the analogues of \accorpa{convtal}{convpal} hold true.
Then, we deduce that \pier{$\btheta_{\alpha_n}\to\btheta$ weakly  in $\C0{\VA\rho}$ 
and $\bphi_{\alpha_n}\to\bphi$ strongly} in $\C0{\Lx r}$ for $1\leq r<6$\juerg{, by virtue of, e.g.,} \cite[Sect.~8, Cor.~4]{Simon},
and it also follows that}
\begin{align}
\label{stark1}
&\gianni{\ell(\bphi_{\alpha_n})\to\ell(\bphi) \aand
\ell'(\bphi_{\alpha_n})\to\ell'(\bphi)} 
\non\\[0.5mm]
&\quad\hbox{\gianni{strongly in $\C0{\Lx r}$ for $1\leq r<6$,}}
\\[0.5mm]
\label{stark3}
&g_1^{\alpha_n}\to \beta_1(\bphi(T)-\phi_\Omega)\quad\mbox{strongly in }\,H,\\[0.5mm]
\label{stark4}
&g_2^{\alpha_n}\to \beta_2(\bphi-\phi_Q)\quad\mbox{strongly in }\,L^2(Q),\\[0.5mm]
\label{stark5}
&g_3^{\alpha_n}\to\beta_3(\btheta(T)-\theta_\Omega)\quad\mbox{weakly in }\,\VA\rho,\\[0.5mm]
\label{stark6}
&g_4^{\alpha_n}\to\beta_4(\btheta-\theta_Q)\quad\mbox{strongly in }\,L^2(Q),
\end{align} 
\gianni{as well as
\Beq
  F_2''(\bphi_{\alpha_n})\to F_2''(\bphi)
  \quad \hbox{strongly in $\LQ r$ for $1\leq r<+\infty$}\,,
  \label{stark7}
\Eeq
since $F_2''$ is continuous and bounded.}
Moreover, by virtue of the estimates \eqref{pqboundal} \pier{and \eqref{lamal},} and invoking
\cite[Sect.~8, Cor.~4]{Simon} once more,
there are limits $\,q,p,\Lambda$ such that, at least for a subsequence which is again
indexed by~$n$,
\begin{align}
\label{limq}
q_{\alpha_n}&\to q\quad\mbox{weakly-star in }\,H^1(0,T;H)\cap L^\infty(0,T;\VA\rho)\cap L^2(0,T;\VA{2\rho}),
\nonumber\\[0.5mm]
&\hspace*{12.8mm}\mbox{weakly in $\,C^0([0,T];\VA\rho)$ and strongly in }\,\pier{C^0}([0,T];H),\\[0.5mm]
\label{limp}
p_{\alpha_n}&\to p\quad \mbox{weakly-star in }\, \pier{L^\infty(0,T;H)\cap L^2(0,T;\Huno)},
\\[0.5mm] 
\label{limlam}
\Lambda^{\alpha_n}&\to\Lambda\quad\mbox{weakly in }\,{\cal Z}^*.
\end{align}
With these convergence results, it is an easy task to show that
\begin{align}\label{Fredi}
&\ell(\bphi_{\alpha_n})p_{\alpha_n}\to\ell(\bphi)p, \quad \ell(\bphi_{\alpha_n})\dt q_{\alpha_n}\to \ell(\bphi)\dt q, \quad
F_2''(\bphi_{\alpha_n})p_{\alpha_n}\to F_2''(\bphi)p,\nonumber\\[0.5mm]
&\ell'(\bphi_{\alpha_n})\btheta_{\alpha_n}p_{\alpha_n}\to\ell'(\bphi)\btheta\, p,\quad\mbox{all weakly in }\,L^1(Q),
\end{align}
\gianni{by using for the latter \eqref{stark1} with $r=4$, 
the strong convergence $\btheta_{\alpha_n}\to\btheta$ in $\C0H$,
and the weak convergence $p_{\alpha_n}\to p$ in $\L\infty{\Lx4}$
ensured by \eqref{convtal} and~\eqref{limp}, respectively.}
A fortiori, since all of the sequences occurring in \eqref{Fredi} are bounded in $L^2(Q)$, we even have weak convergence in~$L^2(Q)$.

At this point, we write the variational inequality \eqref{vugal} for $\alpha=\alpha_n$, $n\in\enne$,
and pass to the limit as $n\to\infty$, which immediately yields the validity of \eqref{vug3}. 
Next, we easily see that the final value \pier{condition \eqref{pier2} holds} true. 
Moreover, writing \eqref{as1al} with 
$\alpha=\alpha_n$, $n\in\enne$, and passing to the limit as $n\to\infty$, we recover~\eqref{as1}.
It  remains to show that \eqref{limas2} is satisfied\pier{, but this can be easily achieved by taking the limit in \eqref{pier1}  written for $\alpha=\alpha_n$, because of \eqref{stark3}--\eqref{stark5} and \eqref{limp}--\eqref{Fredi}. With this, the assertion is proved.}
\Edim 

\Brem 
Unfortunately, we are unable to derive any complementarity slackness 
conditions for the Lagrange multiplier $\Lambda$. Indeed, while it is easily seen that
$$
\liminf_{n\to\infty}\iint_Q \Lambda_{\alpha_n}\,p_{\alpha_n}\,=\,\liminf_{n\to\infty}
\iint_Q \alpha_n \, h''(\bphi_{\alpha_n})\,|p_{\alpha_n}|^2\,\ge\,0
\quad \forall\,n\in\enne,
$$
the available convergence properties do not suffice to conclude that
$\,\langle \Lambda,p\rangle_{\cal Z}\,\ge\,0$.
\Erem


\section*{Acknowledgments}
\pier{This research was supported by the Italian Ministry of Education, 
University and Research~(MIUR): Dipartimenti di Eccellenza Program (2018--2022) 
-- Dept.~of Mathematics ``F.~Casorati'', University of Pavia. 
In addition, PC and CG gratefully acknowledge some other 
financial support from the GNAMPA (Gruppo Nazionale per l'Analisi Matematica, 
la Probabilit\`a e le loro Applicazioni) of INdAM (Isti\-tuto 
Nazionale di Alta Matematica) and the IMATI -- C.N.R. Pavia, Italy.}


\vspace{3truemm}

\Begin{thebibliography}{10}


\bibitem{BCGMR}
V. Barbu, P. Colli, G. Gilardi, G. Marinoschi, E. Rocca, Sliding mode control
for a nonlinear phase-field system, {\em SIAM J. Control Optim.} {\bf 55} (2017), 2108-2133.

\bibitem{Brezis}
H. Brezis,
``Op\'erateurs maximaux monotones et semi-groupes de contractions
dans les espaces de Hilbert'',
North-Holland Math. Stud. Vol.
{\bf 5},
North-Holland,
Amsterdam,
1973.

\bibitem{BS}
M. Brokate, J. Sprekels, ``Hysteresis and Phase Transitions'',
Applied Mathematical Sciences Vol. {\bf 121}, Springer, New York, 1996.

\bibitem{Cp1}
G. Caginalp, An analysis of a phase field model of a free boundary,
{\it Arch. Rat. Mech. Anal.} \textbf{92} (1986), 205--245.

\bibitem{Cp2}
G. Caginalp, Stefan and Hele--Shaw models as asymptotic limits of the phase-field equations,
{\it Phys. Rev. A} {\bf 39} (1989), 5887--5896.

\bibitem{Cp3}
G. Caginalp, The dynamics of a conserved phase field system: Stefan-like, Hele--Shaw, and
Cahn--Hilliard models as asymptotic limits, {\it IMA J. Appl. Math.} {\textbf 44} (1990),
77--94.

\bibitem{Cp4}
G. Caginalp, Phase field models and sharp interfaces: some differences in subtle situations,
{\it Rocky Mountains J. Math.} {\textbf 21} (1991), 603--616.


\bibitem{ChenHoff}
Z. Chen, K.-H. Hoffmann, Numerical solutions of the optimal control problem governed by a phase field model,
in: ``Estimation and Control of Distributed Parameter Systems'' (W. Desch, F. Kappel, K. Kunisch, eds.), 
Int. Ser. Numer. Math. Vol. {\bf 100}, pp. 79--97, Birkh\"auser, Basel,  1991.    

\bibitem{CFS}
P. Colli, M.H. Farshbaf-Shaker, J. Sprekels, A deep quench approach to the optimal  control
of an Allen--Cahn equation with dynamic boundary conditions, {\em Appl. Math. Optim.} {\bf 71}
(2015), 1--24.

\bibitem{CFGS}
P. Colli, M.H. Farshbaf-Shaker, G. Gilardi, J. Sprekels,
Optimal boundary control of a viscous Cahn--Hilliard system
with dynamic boundary condition and double obstacle potentials,
{\it SIAM J. Control Optim.} {\bf 53} (2015), 2696--2721.

\bibitem{CG}
P. Colli, G. Gilardi, Well-posedness, regularity and asymptotic analyses for a fractional phase field system, {\em Asympt. Anal.} {\bf 114} (2019), 93--128. 

\bibitem{CGMR1}
P. Colli, G. Gilardi, G. Marinoschi, E. Rocca,
Optimal control for a phase field system with a possibly singular potential,
{\it Math. Control Relat. Fields\/} {\bf 6} (2016), 95--112.

\bibitem{CGMR2}
P. Colli, G. Gilardi, G. Marinoschi, E. Rocca,
Optimal control for a conserved phase field system with a possibly singular potential,
{\it Evol. Equ. Control Theory\/} {\bf 7} (2018), 95--116.

\bibitem{CGRS}
P. Colli, G. Gilardi, E. Rocca, J. Sprekels,
Optimal distributed control of a diffuse interface model of tumor growth,
\textit{Nonlinearity} \textbf{30} (2017), 2518--2546.




\bibitem{CGS_DQ}
P. Colli, G. Gilardi, J. Sprekels,
Optimal velocity control of a convective Cahn--Hilliard system with double obstacles
and dynamic boundary conditions: a `deep quench' approach,
{\it J. Convex Anal.} {\bf 26} (2019), 485--514.

\bibitem{CGS18}
P. Colli, G. Gilardi, J. Sprekels,
Well-posedness and regularity for a generalized fractional Cahn--Hilliard system,
{\em Atti Accad. Naz. Lincei Rend. Lincei Mat. Appl.} {\bf 30} (2019), 437--478.

\juerg{
\bibitem{CGS21bis}
P. Colli, G. Gilardi, J. Sprekels,
Recent results on well-posedness and optimal control for a class of generalized
fractional Cahn--Hilliard systems, {\em Control Cybernet.} {\bf 48} (2019), 153--197.  
}

\bibitem{CGS19}
P. Colli, G. Gilardi, J. Sprekels,
Optimal distributed control of a generalized fractional Cahn--Hilliard system,
{\em Appl. Math. Optim.}, doi:10.1007/s00245-018-9540-7.

\bibitem{CGS21}
P. Colli, G. Gilardi, J. Sprekels,
Deep quench approximation and optimal control 
of general Cahn--Hilliard systems with fractional operators 
and double obstacle potentials, {\em Discrete Contin. Dyn. Syst. Ser. S},
doi:10.3934/dcdss.2020213.

\bibitem{CGSASY}
P. Colli, G. Gilardi, J. Sprekels,
Asymptotic analysis of a tumor growth model with fractional operators,
{\em Asympt. Anal.}, doi:10.3233/ASY-191578.

\bibitem{CGS25}
P. Colli, G. Gilardi, J. Sprekels,
A distributed control problem for a fractional tumor growth model,
{\em Mathematics} {\bf 2019}, 7, 792.

\pier{\bibitem{CM}
P. Colli, D. Manini, 
Sliding mode control for a generalization 
of the Caginalp phase-field system,
{\it Appl. Math. Optim.}, doi:10.1007/s00245-020-09682-3.}

\bibitem{CMR}
P. Colli, G. Marinoschi, E. Rocca,
Sharp interface control in a Penrose--Fife model, {\em ESAIM: COCV} {\bf 22} (2016), 473--499.

\bibitem{CSS1}
P. Colli, A. Signori, J. Sprekels,
Optimal control of a phase field system modelling tumor growth with chemotaxis and singular
potentials, {\em Appl. Math. Optim.}, doi:10.1007/s00245-019-09618-6.


\bibitem{NPV}
E. Di Nezza, G. Palatucci, E. Valdinoci, 
Hitchhiker's guide to the fractional Sobolev spaces, {\em  Bull. Sci. Math.} {\bf 136} (2012), 521--573.

\bibitem{EK}
M. Ebenbeck, P. Knopf,
Optimal medication for tumors modeled by a Cahn--Hilliard--Brinkman equation,
\pier{{\em Calc. Var. Partial Differential Equations}  {\bf 58}  (2019),  Paper No. 131, 31 pp.}

\bibitem{EK_ADV}
M. Ebenbeck, P. Knopf,
Optimal control theory and advanced optimality conditions 
for a diffuse interface model of tumor growth,
preprint arXiv:1903.00333 [math.OC] (2019), 1--34.

\bibitem{GLR}
H. Garcke, K. F. Lam, E. Rocca,
Optimal control of treatment time in a diffuse interface model for tumour growth,
\textit{Appl. Math. Optim.} {\bf 78}, 495--544.

\bibitem{Hein}
M. Heinkenschloss, The numerical solution of a control problem governed by a phase field model,
{\em Optim. Methods Softw.} {\bf 7} (1997), 211--263.

\bibitem{HeinSachs}
M. Heinkenschloss, E. W. Sachs, Numerical solution of a constraint control problem for a phase field model,
in: ``Control and Estimation of Distributed Parameter Systems'' (W. Desch, F. Kappel, K. Kunisch, eds.), 
Int. Ser. Numer. Math. Vol. {\bf 118}, pp. 171--188, Birkh\"auser, Basel,  1994.    

\bibitem{HeinTr}
M. Heinkenschloss, F. Tr\"oltzsch, Analysis of the Lagrange--SQP--Newton method for the control of a phase field
equation, {\em Control Cybernet.} {\bf 28} (1999), 178--211.


\bibitem{HoffJiang}
K.-H. Hoffmann, L. Jiang, 
Optimal control problem of a phase field model for solidification, {\em Numer. Funct. Anal. Optim.} {\bf 13} (1992),
11--27. 

\bibitem{HSS}
W. Horn, J. Sokolowski, J. Sprekels, 
A control problem with state constraints for the Penrose--Fife phase-field model, {\em Control Cybern.} {\bf 25} (1996), 
1137--1153. 
  
\bibitem{LSU}
O.A. Lady\v zenskaja, V.A. Solonnikov and N.N. Uralceva,
``Linear and quasilinear equations of parabolic type'',
{\it Mathematical Monographs Volume} {\bf 23},
American Mathematical Society, Providence, Rhode Island, 1968.

\bibitem{LS}
C. Lefter, J. Sprekels, Control of a phase field system modeling non-isothermal phase transitions,
{\em Adv. Math. Sci. Appl.} {\bf 17} (2007), 181--194.


\bibitem{S_a}
A. Signori,
Vanishing parameter for an optimal control problem modeling tumor growth,
{\em Asympt. Anal.} \pier{{\bf 117} (2020), 43--66}.                        

\bibitem{S_DQ}
A. Signori,
Optimality conditions for an extended tumor growth model with 
double obstacle potential via deep quench approach, 
\pier{{\em Evol. Equ. Control Theory} {\bf 9} (2020), 193--217.}

\bibitem{S}
A. Signori,
Optimal distributed control of an extended model of tumor 
growth with logarithmic potential,
{\it Appl. Math. Optim.}, doi:10.1007/s00245-018-9538-1.

\bibitem{S_b}
A. Signori,
Optimal treatment for a phase field system of Cahn--Hilliard 
type modeling tumor growth by asymptotic scheme, 
{\em Math. Control Rel. Fields}, doi:10.3934/mcrf.2019040.

\bibitem{Simon}
J. Simon,
{Compact sets in the space $L^p(0,T; B)$},
{\it Ann. Mat. Pura Appl.~(4)\/} 
{\bf 146} (1987) 65--96.


\bibitem{SZheng}
 J. Sprekels, S. Zheng, Optimal control problems for a thermodynamically consistent 
  model  of  phase-field type for  phase  transitions, {\em  Adv.  Math.  Sci.  Appl.} {\bf 1} (1992), 
	113--125.

\bibitem{Tr}
F. Tr\"oltzsch, ``Optimal Control of Partial Differential Equations: Theory, Methods and Applications'',
Graduate Studies in Mathematics Vol. {\bf 12}, American Mathematical Society, Providence, Rhode Island, 2010.

\End{thebibliography}

\End{document}
